\documentclass[12pt,reqno]{amsart}
\usepackage{amsmath,amssymb,amsthm}
\textwidth=14cm\textheight=23cm \numberwithin{equation}{section}
\theoremstyle{plain}
\newtheorem{theorem}{Theorem}[section]

\newtheorem{proposition}{Proposition}[section]
\newtheorem{lemma}{Lemma}[section]
\newtheorem{corollary}{Corollary}[section]
\theoremstyle{definition}
\newtheorem{definition}{Definition}[section]

\newtheorem*{prob1}{Problem}
\theoremstyle{remark}

\newtheorem*{rem0}{Remark}

\newtheorem*{acknow}{Acknowledgments}

\begin{document}
\title[Lagrangian submanifolds in complex hyperquadrics]
{On Lagrangian submanifolds \\in complex hyperquadrics and \\
isoparametric hypersurfaces in spheres}
\author{Hui Ma}
\address{Department of Mathematical Sciences, Tsinghua University,
Beijing 100084, P.R. CHINA} \email{hma@math.tsinghua.edu.cn}
\thanks{The first named author is partially supported by NSFC grant
No.~10501028, SRF for ROCS, SEM and NKBRPC No.~2006CB805905. The
second named author is partially supported by JSPS grant of Kiban
(A) No.~17204006}
\author{Yoshihiro Ohnita}
\address{Department of Mathematics, Osaka City University, Sugimoto,
Sumiyoshi-ku, Osaka, 558-8585, JAPAN}
\email{ohnita@sci.osaka-cu.ac.jp} \dedicatory{Dedicated to Professor
Hajime Urakawa on his sixtieth birthday}
\date{August 10, 2008}

\begin{abstract}
The $n$-dimensional complex hyperquadric is a compact complex
algebraic hypersurface defined by the quadratic equation in the
$(n+1)$-dimensional complex projective space, which is isometric to
the real Grassmann manifold of oriented $2$- planes and is a compact
Hermitian symmetric space of rank $2$. In this paper we study
geometry of compact Lagrangian submanifolds in complex hyperquadrics
from the viewpoint of the theory of isoparametric hypersurfaces in
spheres. From this viewpoint we provide a classification theorem of
compact homogeneous Lagrangian submanifolds in complex hyperquadrics
by using the moment map technique. Moreover we determine the
Hamiltonian stability of compact minimal Lagrangian submanifolds
embedded in complex hyperquadrics which are obtained as Gauss images
of isoparametric hypersurfaces in spheres with $g(=1,2,3)$ distinct
principal curvatures.
\end{abstract}

\maketitle

\section*{Introduction}
\label{intro}

Let $(M,\omega)$ be a $2n$-dimensional symplectic manifold with
symplectic form $\omega$. The {\it Lagrangian submanifold} $L$ in
$M$ is by definition an $n$-dimensional smooth submanifold $L$ in
$M$ such that the pull-back of $\omega$ to $L$ vanishes. It is an
interesting problem to investigate Lagrangian submanifolds in
specific symplectic manifolds such as K\"ahler manifolds, complex
space forms, Hermitian symmetric spaces, generalized flag manifolds
with invariant symplectic forms, toric K\"ahler manifolds etc., from
both viewpoints of symplectic geometry and Riemannian geometry
(cf.~Section \ref{sec:1}).

In this paper we study compact Lagrangian submanifolds in complex
hyperquadrics $Q_{n}({\mathbf C})$. The $n$-dimensional complex
hyperquadric $Q_{n}({\mathbf C})$ is a compact complex algebraic
hypersurface defined by the quadratic equation
$z_{0}^{2}+z_{1}^{2}+\cdots+z_{n}^{2}+z_{n+1}^{2}=0$ in the
$(n+1)$-dimensional complex projective space, which is isometric to
the real Grassmann manifold of oriented $2$-dimensional vector
subspaces of ${\mathbf R}^{n+2}$ and is a compact Hermitian
symmetric space of rank $2$. In Section~\ref{sec:2}, we discuss a
link of Lagrangian Geometry in the complex hyperquadrics with
Hypersurface Geometry in the sphere. We know a fundamental fact that
the Gauss map of any oriented hypersurface in the standard sphere
$S^{n+1}(1)$ is always a Lagrangian immersion into the complex
hyperquadric $Q_{n}({\mathbf C})$. It follows from the mean
curvature formula of B.~Palmer (\cite{Palmer97}) that the images of
the Gauss map of compact oriented hypersurfaces with constant
principal curvatures, so called {\it isoparametric hypersurfaces},
in spheres provide a nice class of compact minimal Lagrangian
submanifolds embedded in complex hyperquadrics. Particularly
homogeneous isoparametric hypersurfaces in spheres correspond to
minimal Lagrangian orbits in complex hyperquadrics.
The theory of isoparametric hypersurfaces was originated by \'{E}lie
Cartan.
We know that by the result of W.~Y.~Hsiang-J.~B.~Lawson
(\cite{Hsiang-Lawson1971}) homogeneous isoparametric hypersurface
can be obtained as a principal orbit of the isotropy representation
of a compact Riemannian symmetric pair $(U,K)$ of rank $2$. Since
non-homogeneous isoparametric hypersurfaces in spheres exist
(H.~Ozeki-M.~Takeuchi (\cite{Ozeki-TakeuchiI,Ozeki-TakeuchiII})) and
are almost classified (D.~Ferus-H.~Karcher-H.~F.~M\"unzner
(\cite{FKM}), T.~Cecil-Q.~S.~Chi-G.~R.~Jensen
(\cite{Cecil-Chi-Jensen})) at present, they also give many nice
examples of compact non-homogeneous minimal Lagrangian submanifolds
in complex hyperquadrics.

As the theory of isoparametric hypersurfaces are now well-developed
in differential geometry, we are interested in the following problem
:

\begin{prob1}
Investigate properties of compact minimal Lagrangian submanifolds in
complex hyperquadrics obtained as Gauss images of isoparametric
hypersurfaces in spheres. More generally, study compact Lagrangian
submanifolds in complex hyperquadrics by using the theory of
isoparametric hypersurfaces in spheres.
\end{prob1}

A Lagrangian submanifold obtained as a Lagrangian orbit of the
Hamiltonian group action on $M$ is called a {\it homogeneous}
Lagrangian submanifold. The investigation of the Lie-theoretic
construction of Lagrangian submanifolds is inspired by
\cite{Oh94},\cite{BedulliGori06-1}. In Section \ref{sec:3} we shall
discuss Lie algebraic properties of homogeneous isoparametric
hypersurfaces in spheres related with homogeneous Lagrangian
geometry in complex hyperquadrics such as the moment maps.

In Section \ref{sec:4}, based on the link with theory of homogeneous
isoparametric hypersurfaces, we shall show a classification theorem
of compact homogeneous Lagrangian submanifolds, i.e.,~Lagrangian
orbits of compact connected Lie subgroups of $SO(n+2)$, in
$Q_{n}({\mathbf C})$. The moment map technique plays an essential
role in the argument. We shall prove that all compact homogeneous
Lagrangian submanifolds in $Q_{n}({\mathbf C})$ are the Gauss image
of compact homogeneous isoparametric hypersurfaces in spheres, or
are obtained as their Lagrangian deformations in the following cases
: Suppose that $(U,K)$ is one of
\begin{enumerate}
\renewcommand{\labelenumi}{(\alph{enumi})}
\item $(S^{1}\times SO(3),SO(2))$,
\item $(SO(3)\times SO(3),SO(2)\times SO(2))$,
\item $(SO(3)\times SO(n+1),SO(2)\times SO(n)){\ }(n\geq{3})$,
\item $(SO(m+2),SO(2)\times SO(m)){\ }(n=2m-2,m\geq{3})$.
\end{enumerate}
In the first two cases, it is elementary and well-known to describe
all Lagrangian orbits of the natural actions of $K=SO(2)$ on
$Q_{1}({\mathbf C})\cong S^{2}$ and $K=SO(2)\times SO(2)$ on
$Q_{2}({\mathbf C})\cong{S^{2}\times S^{2}}$. Also in the last two
cases there exist one-parameter families of Lagrangian $K$-orbits in
$Q_{n}({\mathbf C})$ and each family contains Lagrangian
submanifolds which can NOT be obtained as the Gauss image of any
homogeneous isoparametric hypersurface in a sphere. The fourth one
is a new family of Lagrangian orbits and it will be discussed in detail.

The notion of Hamiltonian minimality and Hamiltonian stability for
Lagrangian submanifolds in K\"ahler manifolds was introduced and
investigated first by Y.~G.~Oh
(\cite{Oh90},\cite{Oh91},\cite{Oh93},\cite{Oh94}). A compact
Lagrangian submanifold $L$ in a K\"ahler manifold $M$ is called {\it
Hamiltonian minimal}{\ }(shortly, {\it H-minimal}) if it has
extremal volume under every Hamiltonian deformations of $L$. A
compact H-minimal Lagrangian submanifold in a K\"ahler manifold $M$
is called {\it Hamiltonian stable}{\ }(shortly, {\it H-stable}) if
the second variation for the volume is nonnegative under every
Hamiltonian deformations of $L$.
We do NOT know so many examples of compact Hamiltonian stable
Lagrangian submanifolds yet.

\begin{prob1}
Construct and classify compact Hamiltonian stable minimal or
H-minimal Lagrangian submanifolds in specific K\"ahler manifolds.
\end{prob1}

In Section \ref{sec:5} we shall determine Hamiltonian stability of
compact minimal Lagrangian submanifolds embedded in $Q_{n}({\mathbf
C})$ obtained as Gauss images of isoparametric hypersurfaces with
$g=1,2,3$. Due to the result of E.~Cartan, they all are homogeneous.
The Gauss images of compact isoparametric hypersurfaces with $g=2$
are not Hamiltonian stable if and only if the multiplicities
$m_{1}$, $m_{2}$ of the principal curvatures satisfy
$m_{2}-m_{1}\geq{3}$. We shall show that the Gauss images of all
compact isoparametric hypersurfaces with $g=3$ are Hamiltonian
stable and they provide new examples of compact Hamiltonian stable
minimal Lagrangian submanifolds embedded in $Q_{n}({\mathbf C})$.

\begin{acknow}
This work has been done during the second named author's visits at
Tsinghua University in Beijing and Fudan University in Shanghai in
2005-2006 and the first named author's stay at Osaka City University
Advanced Mathematical Institute (OCAMI) in 2007. The authors would
like to thank Professors Zizhou Tang and Yuanlong Xin for helpful
conversations and generous supports, and Professor Akio Kawauchi
(Director of OCAMI) for excellent support and research circumstance
based on OCAMI. They also would like to thank Professors Franz Pedit
and Ulrich Pinkall for suggestion on the elementary version of
Theorem \ref{HamiltonianDeformationsIsomonodromyDeformation}.
\end{acknow}


\section{Hamiltonian deformations of Lagrangian submanifolds}
\label{sec:1}

Let $(M,\omega)$ be a $2n$-dimensional symplectic manifold with a
symplectic form $\omega$. A {\it Lagrangian immersion}
$\varphi:L\longrightarrow M$ is a smooth immersion of an
$n$-dimensional smooth manifold $L$ into $M$ satisfying the
condition $\varphi^{\ast}\omega=0$. Generally a smooth immersion
$\varphi:Q\longrightarrow M$ of a $p$-dimensional smooth manifold
$Q$ into $M$ satisfying the condition $\varphi^{\ast}\omega=0$ is
called {\it isotropic}.

The normal bundle $\varphi^{-1}TM/\varphi_{\ast}TL$ of a Lagrangian
immersion $\varphi:L\longrightarrow M$ can be identified with the
cotangent bundle $T^{\ast}L$ of $L$ :
$$
\varphi^{-1}TM/\varphi_{\ast}TL \ni v\longmapsto
\alpha_{v}:=\omega(v,\cdot)\in T^{\ast}L
$$
By the definition a {\it Lagrangian deformation} is a smooth family
of Lagrangian immersions $\varphi_{t}:L\longrightarrow M$ with
$\varphi=\varphi_{0}$. Let $\alpha_{V_{t}}$ be a $1$-form on $L$
corresponding to its variational vector field
$$
V_{t}:=\frac{\partial{\varphi_{t}}}{\partial{t}}\in
C^{\infty}(\varphi_{t}^{-1}TM).
$$
The Lagrangian deformation is characterized by the condition that
$\alpha_{V_{t}}$ is closed, i.e., $\alpha_{V_{t}}\in Z^{1}(L)$, for
each $t$. Furthermore, if $\alpha_{V_{t}}$ is exact, i.e.,
$\alpha_{V_{t}}\in B^{1}(L)$, for each $t$, then $\{\varphi_{t}\}$
is called a {\it Hamiltonian deformation} of $\varphi=\varphi_{0}$.

Suppose that $[(1/2\pi)\omega]\in H^{2}(M;{\mathbf R})$ is an
integral cohomology class. Then there is a complex line bundle
${\mathcal L}$ over $M$ with a $U(1)$-connection $\nabla$ whose
curvature coincides with $\sqrt{-1}\omega$, and thus the Lagrangian
property is equivalent to the flatness of the induced connection in
the pull-back $U(1)$-bundle $\varphi^{-1}{\mathcal L}$. Then we can
show

\begin{theorem}\label{HamiltonianDeformationsIsomonodromyDeformation}
A Lagrangian deformation $\{\varphi_{t}\}$ is a Hamiltonian
deformation if and only if $\{\varphi_{t}\}$ provides an {\it
isomonodromy deformation} of the induced flat $U(1)$-connections in
$\varphi_{t}^{-1}{\mathcal L}$, that is, they have the same holonomy
homomorphism $\pi_{1}(L)\rightarrow U(1)$.
\end{theorem}

\begin{proof}
For each $t$ let $\rho_{t}:\pi_{1}(L)\rightarrow{U(1)}$ denote the
holonomy homomorphism of the induced flat $U(1)$-connection
$\nabla^{t}=\varphi_{t}^{-1}\nabla$ in $\varphi_{t}^{-1}{\mathcal
L}$. Then by straightforward calculations we obtain the formula
\begin{equation}
\rho_{t}([\gamma])^{-1}\frac{d}{dt}\rho_{t}([\gamma])
=-\sqrt{-1}\int_{\gamma}\alpha_{V_{t}}
\end{equation}
for each $[\gamma]\in{\pi_{1}(L)}$. Hence the statement of the
theorem follows from this formula.
\end{proof}

\subsection{Moment maps and Lagrangian orbits}
We call a Lagrangian submanifold obtained as a Lagrangian orbit of a
Lie group {\it a homogeneous Lagrangian submanifold}.
Suppose that a connected Lie group $K$ acts on $M$ in a Hamiltonian
way. Let $\mathfrak{k}$ denote the Lie algebra of Lie group $K$ and
${\mathfrak k}^{\ast}$ its dual vector space. Let $\mu :
M\longrightarrow{\mathfrak k}^{\ast}$ be the moment map for the
Hamiltonian group action $K$ on $M$. By the definition the moment
map $\mu$ satisfies the following conditions :
\begin{enumerate}
\renewcommand{\labelenumi}{(\arabic{enumi})}
\item
$d\langle{\mu},{\xi}\rangle=\omega(\tilde{\xi},\cdot)$ for all
$\xi\in{\mathfrak k}$.
\item
$\mu(a\cdot{x})=\mathrm{Ad}^{\ast}(a^{-1})\mu(x)$ for all $x\in M$
and all $a\in K$.
\end{enumerate}
Here $\tilde{\xi}$ denotes the vector field on $M$ induced by the
action of the one-parameter subgroup $\mathrm{exp}(t\xi)$. Set
$$
{\mathfrak z}({\mathfrak k}^{\ast}):=\{\alpha\in{\mathfrak
k}^{\ast}{\ }\vert{\ } \mathrm{Ad}^{\ast}(a)\alpha=\alpha \mbox{ for
all } a\in K\}.
$$
Then we know the following basic properties on moment maps and
Lagrangian orbits (see also \cite{Biliotti06}).

\begin{proposition}
Let $x\in{M}$. If $K\cdot{x}$ is an isotropic orbit, then the orbit
$K\cdot{x}$ is contained in a level set $\mu^{-1}(\alpha)$ of the
moment map $\mu$ for some $\alpha\in{\mathfrak k}^{\ast}$, and
$\alpha\in{\mathfrak z}({\mathfrak k}^{\ast})$.
\end{proposition}

\begin{proposition}
For each $\alpha\in{\mathfrak z}({\mathfrak k}^{\ast})$, choose an
arbitrary point $x\in\mu^{-1}(\alpha)$. Then the $K$-orbit
$K\cdot{x}$ in $M$ through $x$ has the following properties :
\begin{enumerate}
\renewcommand{\labelenumi}{(\arabic{enumi})}
\item
The orbit $K\cdot{x}$ of $K$ through $x$ is contained in
$\mu^{-1}(\alpha)$.
\item
The orbit $K\cdot{x}$ is an isotropic submanifold of $M$.
\item
$K\cdot{x}$ is a Lagrangian submanifold of $M$ if and only if
$$
T_{y}(K\cdot{x})=\mathrm{Ker}(d\mu)_{y}\quad \mbox{ for each } y\in
K\cdot{x},
$$
that is, $K\cdot{x}$ is an open subset of $\mu^{-1}(\alpha)$.
\item
Assume that the action of $K$ on $M$ is proper. Then $K\cdot{x}$ is
a Lagrangian submanifold of $M$ if and only if the orbit $K\cdot{x}$
is a connected component of $\mu^{-1}(\alpha)$.
\end{enumerate}
\end{proposition}

In the case when $K$ is compact, we can identified ${\mathfrak
k}^{\ast}$ with ${\mathfrak k}$ via a fixed
$\mathrm{Ad}(K)$-invariant inner product of ${\mathfrak k}$.
The following results are on the connectivity of the level subset
and the convexity of the image of the moment maps.

\begin{proposition}[\cite{Kirwan84}, cf.~\cite{Guillemin-Sjamaar05}]
\label{KirwanResult} Suppose that $K$ is a compact, connected Lie
group with the Hamiltonian action on a compact symplectic manifold
$M$. Let $T$ be a maximal torus of $K$ with Lie algebra ${\mathfrak
t}$.
\begin{enumerate}
\renewcommand{\labelenumi}{(\arabic{enumi})}
\item
For each $\xi\in{\mathfrak c}({\mathfrak k})\cong{\mathfrak
z}({\mathfrak k}^{\ast})$, the level set $\mu^{-1}(\xi)$of the
moment map $\mu$ is a connected subset of $M$.
\item
The intersection $\mu(M)\cap{\mathfrak t}_{+}$ of the image of the
moment map with the positive Weyl chamber ${\mathfrak t}_{+}$ is a
compact connected convex subset of ${\mathfrak k}\cong{\mathfrak
k}^{\ast}$.
\end{enumerate}
\end{proposition}

\begin{prob1}
Classify homogeneous Lagrangian submanifolds, i.e., Lagrangian
orbits of Hamiltonian group actions, in specific symplectic
manifolds.
\end{prob1}


\subsection{Hamiltonian minimality and Hamiltonian stability}

Now we assume that $(M,\omega,J,g)$ is a K\"ahler manifold with
complex structure $J$ and K\"ahler metric $g$. Let
$\varphi:L\rightarrow M$ be a Lagrangian immersion into $M$. Let $H$
denote the mean curvature vector field of $\varphi$ and we call the
corresponding $1$-form $\alpha_{H}\in\Omega^{1}(L)$ the {\it mean
curvature form} of $\varphi$.

The mean curvature form $\alpha_{H}$ must satisfy the following
identity (\cite{Dazord}), which follows from the Codazzi equation :
$d\alpha_{H}=\varphi^{\ast}\rho_{M}$, where $\rho_{M}$ denotes the
Ricci form of $M$. Thus if $M$ is an Einstein-K\"ahler manifold,
then we have $d\alpha_{H}=0$.

Here we suppose that $L$ is compact without boundary.

\begin{definition}
A Lagrangian immersion $\varphi$ is called {\it Hamiltonian minimal}
(shortly, H-minimal) or {\it Hamiltonian stationary} if under every
Hamiltonian deformation $\{\varphi_{t}\}$ the first variation of the
volume vanishes.
\end{definition}

By the first variational formula, the H-minimality equation is
$\delta\alpha_{H}=0$, where $\delta$ denotes the codifferential
operator with respect to the induced metric $\varphi^{\ast}g$ on
$L$. Thus if $M$ is an Einstein-K\"ahler manifold, then the mean
curvature form $\alpha_{H}$ is a harmonic $1$-form on $L$.

\begin{definition}
An H-minimal Lagrangian immersion $\varphi$ is called {\it
Hamiltonian stable} (shortly, {\it H-stable}) if under every
Hamiltonian deformation $\{\varphi_{t}\}$ the second variation of
the volume is nonnegative.
\end{definition}
The second variational formula is as follows (\cite{Oh93}) :
\begin{equation*}
\begin{split}
&\frac{d^{2}}{dt^{2}}\mathrm{Vol}(L,\varphi_{t}^{\ast}g)\vert_{t=0}\\
=&\int_{L}\left( \langle\Delta^{1}_{L}\alpha,\alpha\rangle
-\langle\bar{R}\alpha,\alpha\rangle
-2\langle\alpha\otimes\alpha\otimes\alpha_{H},S\rangle
+\langle\alpha_{H},\alpha\rangle^{2} \right)dv,
\end{split}
\end{equation*}
where we set $\alpha=\alpha_{V_{0}}\in B_{1}(L)$. Here
$$
\langle\bar{R}\alpha,\alpha\rangle
:=\sum^{n}_{i,j=1}\mathrm{Ric}^{M}(e_{i},e_{j})\alpha(e_{i})\alpha(e_{j}),
$$
where $\{e_{i}\}$ is a local orthonormal frame on $L$  and
$$
S(X,Y,Z):=g(JB(X,Y),Z)=\omega(B(X,Y),Z)
$$
for each $X,Y,Z\in{TL}$, which is a symmetric $3$-tensor field on
$L$ defined by the second fundamental form $B$ of $L$ in $M$.

We know that if a K\"ahler manifold $M$ is simply connected or
$b_{1}(M)=0$, then any holomorphic Killing vector field $V$ on $M$
is always a Hamiltonian vector field, and thus it generates a
volume-preserving Hamiltonian deformation of every H-minimal
Lagrangian immersion into $M$. Such a Hamiltonian deformation should
be considered as a trivial element of the null-space of the second
variations.

\begin{definition}
An H-minimal Lagrangian immersion $\varphi$ is called {\it strictly
Hamiltonian stable} (shortly, {\it strictly H-stable}) if $\varphi$
is Hamiltonian stable and the nullspace of the second variations is
exactly the span of the normal projections of holomorphic Killing
vector fields of the ambient K\"ahler manifold $M$.
\end{definition}

In the case when $L$ is a compact minimal Lagrangian submanifold
(i.~e.~$H\equiv{0}$) in an Einstein-K\"ahler manifold $M$ with
Einstein constant $\kappa$, the second variational formula is quite
simplified and it follows that $L$ is H-stable if and only if the
first (positive) eigenvalue $\lambda_{1}$ of the Laplacian
$\Delta^{0}_{L}$ of $L$ on functions satisfies the inequality
$\lambda_{1}\geq \kappa$.

We do NOT know so many examples of compact Hamiltonian stable Lagrangian
submanifolds yet.
The elementary examples of compact
H-stable minimal or H-minimal Lagrangian submanifolds are as follows
: (1) circles on a plane $S^{1}\subset{\mathbf C}$, (2) great
circles and small circles $S^{1}\subset S^{2}={\mathbf C}P^{1}$, (3)
closed circles $S^{1}\subset H^{2}={\mathbf C}H^{1}$, (4) real
projective subspaces ${\mathbf R}P^{n}\subset{\mathbf C}P^{n}$
(\cite{Oh90}), (5) a product of $n+1$ circles
$S^1(r_{0})\times\cdots\times S^{1}(r_{n})\subset {\mathbf C}^{n+1}$
and the quotient space by the $S^1$-action $T^{n}\subset {\mathbf
C}P^{n}$ (\cite{Oh93}).

In \cite{Amar-Ohn03}, the Hamiltonian stability of compact
irreducible minimal Lagrangian submanifolds with parallel second
fundamental form (i.e., $\nabla{S}=0$) embedded in complex
projective spaces was shown : (a) $SU(p)/SO(p){\mathbf
Z}_{p}\subset{\mathbf C}P^{(p-1)(p+2)/2}$, (b) $SU(p)/{\mathbf
Z}_{p}\subset{\mathbf C}P^{p^{2}-1}$, (c) $SU(2p)/Sp(p){\mathbf
Z}_{2p}\subset{\mathbf C}P^{(p-1)(2p+1)}$, (d) $E_{6}/F_{4}{\mathbf
Z}_{3}\subset{\mathbf C}P^{26}$.

Lagrangian submanifolds satisfying $\nabla{S}=0$ are called {\it
parallel Lagrangian submanifolds}. Parallel Lagrangian submanifolds
in complex space forms were classified by H.~Naitoh and M.~Takeuchi
(\cite{Naitoh81,Naitoh-Takeuchi82,Naitoh83I,Naitoh83II}). Recently
the above result of \cite{Amar-Ohn03} is generalized as follows
(\cite{Amar-Ohn02pp,Amar-Ohn06pp}): If $L$ is a compact parallel
Lagrangian submanifold embedded in a complex space form (=${\mathbf
C}P^{n}$, ${\mathbf C}^{n}$ or ${\mathbf C}H^{n}$), then $L$ is
Hamiltonian stable.

More recently, an example of a compact Hamiltonian stable minimal
Lagrangian submanifold in ${\mathbf C}P^{3}$ with $\nabla{S}\not=0$,
which is obtained as a minimal Lagrangian $SU(2)$-orbit in ${\mathbf
C}P^{3}$, was shown by L.~Bedulli and A.~Gori
(\cite{BedulliGori06-2}), and independently by \cite{Ohnita06}: (e)
The orbit $\rho_{3}(SU(2))[z_{0}^{3}+z_{1}^{3}]\subset{\mathbf
C}P^{3}$ by the irreducible unitary representation of $SU(2)$ of
degree $3$ (cf. Section 5) is a $3$-dimensional compact embedded
Hamiltonian stable minimal Lagrangian submanifold with
$\nabla{S}\not=0$.

L.~Bedulli and A.~Gori (\cite{BedulliGori06-1}) characterized the
existence of Lagrangian orbits in compact K\"ahler manifolds with
$b^{1,1}(M)=1$ in terms of the Stein property of their complexified
orbits. Applying the classification theory of \lq\lq prehomogeneous
vector spaces\rq\rq due to M.~Sato and T.~Kimura \cite{SatoKimura},
they classified compact homogeneous Lagrangian submanifolds in
${\mathbf C}P^{n}$ obtained as Lagrangian orbits of compact simple
Lie subgroups of $SU(n+1)$ :
$ 16 \mbox{ examples } ={\ }
[5 \mbox{ examples with }\nabla{S}=0:{\mathbf R}P^{n},
\mathrm{(a)}\sim\mathrm{(d)}]
+
[11 \mbox{ examples with}\nabla{S}\not=0 \ni \mathrm{(e)}] $.

M.~Takeuchi (\cite{Takeuchi84}) classified all compact totally
geodesic Lagrangian submanifolds in compact irreducible Hermitian
symmetric spaces. He proved that they all are real forms of
Hermitian symmetric spaces, i.e., the fixed point subset of
anti-holomorphic isometries, and are given as symmetric R-spaces
canonically embedded in compact Hermitian symmetric spaces. The
Hamiltonian stability of all compact totally geodesic Lagrangian
submanifolds embedded in compact irreducible Hermitian symmetric
spaces (with Einstein constant $1/2$) are known as follows
(\cite{Takeuchi84},\cite{Amar-Ohn03}) :
\smallskip

\begin{center}
Table $1$
\end{center}
\begin{center}
\begin{tabular}
{|c|c|c|c|c|c|} \hline $M$ &$L$ &Einstein &${\lambda}_{1}$ &H-stable
&stable
\\
\hline $G_{p,q}({\mathbf C}),p\leq q$ &$G_{p,q}({\mathbf R})$ &Yes
&$\frac{1}{2}$ &Yes &No
\\
\hline $G_{2p,2q}({\mathbf C}),p\leq q$ &$G_{p,q}({\mathbf H})$ &Yes
&$\frac{1}{2}$ &Yes &Yes
\\
\hline $G_{m,m}({\mathbf C})$ &$U(m)$ &No &$\frac{1}{2}$ &Yes &No
\\
\hline $\displaystyle{SO(2m)}/{U(m)}$ &$SO(m),m\geq 5$ &Yes
&$\frac{1}{2}$ &Yes &No
\\
\hline $\displaystyle{SO(4m)}/{U(2m)},m\geq 3$
&$\displaystyle{U(2m)}/{Sp(m)}$ &No &$\frac{m}{4m-2}$ &No &No
\\
\hline $\displaystyle{Sp(2m)}/{U(2m)}$ &$Sp(m),m\geq 2$ &Yes
&$\frac{1}{2}$ &Yes &Yes
\\
\hline $\displaystyle{Sp(m)}/{U(m)}$ &$\displaystyle{U(m)}/{O(m)}$
&No &$\frac{1}{2}$ &Yes &No
\\
\hline $Q_{p+q-2}({\mathbf C}), q-p\geq 3$ &$Q_{p,q}({\mathbf R})$,
$p\geq 2$ &No &$\frac{p}{p+q-2}$ &No &No
\\
\hline $Q_{p+q-2}({\mathbf C}), 0\leq{q-p}<{3}$ &$Q_{p,q}({\mathbf
R}), p\geq 2$ &No &$\frac{1}{2}$ &Yes &No
\\
\hline $Q_{q-1}({\mathbf C})$, $q\geq 3$ &$Q_{1,q}({\mathbf R})$
&Yes &$\frac{1}{2}$ &Yes &Yes
\\
\hline $\displaystyle{E_{6}}/{T\cdot Spin(10)}$ &$P_{2}({\mathbf
K})$ &Yes &$\frac{1}{2}$ &Yes &Yes
\\
\hline $\displaystyle{E_{6}}/{T \cdot Spin(10)}$
&$\displaystyle{G_{2,2}({\mathbf H})}/{{\mathbf Z}_{2}}$ &Yes
&$\frac{1}{2}$ &Yes &No
\\
\hline $\displaystyle{E_{7}}/{T \cdot E_{6}}$
&$\displaystyle{SU(8)}/{Sp(4){\mathbf Z}_{2}}$ &Yes &$\frac{1}{2}$
&Yes &No
\\
\hline $\displaystyle{E_{7}}/{T\cdot E_{6}}$ &$\displaystyle{T\cdot
E_{6}}/{F_{4}}$ &No &$\frac{1}{6}$ &No &No
\\
\hline
\end{tabular}
\end{center}
where $G_{p,q}(\mathbf{F})$ : Grassmanian manifold of all
$p$-dimensional subspaces of $\mathbf{F}^{p+q}$, for each
$\mathbf{F}=\mathbf{R}, \mathbf{C}, \mathbf{H}$. $P_2(\mathbf{K})$ :
Cayley projective plane. $Q_n(\mathbf{C})$ : complex hyperquadric of
complex dimension $n$. Here each $M$ is equipped with the standard
K\"ahler metric of Einstein constant $1/2$ and $\lambda_{1}$ denotes
the first eigenvalue of the Laplacian of $L$ on smooth functions.

\begin{rem0}
The second named author apologizes that there are some inaccuracies in
the cases of $M=Q_n(\mathbf{C})$ at the table of
\cite[p.608]{Amar-Ohn03}.
It should be corrected as above.
\end{rem0}

\section{Lagrangian submanifolds in complex hyperquadrics
and hypersurface geometry in spheres} \label{sec:2}

Next we shall discuss Lagrangian submanifolds in the complex
hyperquadrics
$$ Q_{n}({\mathbf C}) \cong
\widetilde{\mathrm{Gr}}_{2}({\mathbf R}^{n+2}) \cong
SO(n+2)/SO(2)\times SO(n),
$$
the latter are compact irreducible Hermitian symmetric spaces of
rank $2$ if $n\geq{3}$ and $S^{2}\times S^{2}$ if $n=2$. Here
$Q_{n}({\mathbf C})$ denotes the complex hypersurface of ${\mathbf
C}P^{n+1}$ defined by the algebraic equation
$z_{0}^{2}+z_{1}^{2}+\cdots+z_{n+1}^{2}=0$ and
$\widetilde{\mathrm{Gr}}_{2}({\mathbf R}^{n+2})$ denotes the real
Grassmann manifold of oriented $2$-planes in ${\mathbf R}^{n+2}$.
Let ${\mathcal L}$ be the tautological holomorphic line bundle over
$Q_{n}({\mathbf C})$ and through the identification $Q_{n}({\mathbf
C})\cong\widetilde{\mathrm{Gr}}_{2}({\mathbf R}^{n+2})$, the
holomorphic line bundle ${\mathcal L}$ can be also considered as the
tautological real vector bundle $\mathcal{V}$ of rank $2$ over
$\widetilde{\mathrm{Gr}}_{2}({\mathbf R}^{n+2})$.

Geometry of Lagrangian submanifolds in complex hyperquadrics has the
important relationship with {\it Hypersurface Geometry} in the unit
sphere $S^{n+1}(1)$. Let $N^{n}\subset S^{n+1}(1)\subset{\mathbf
R}^{n+2}$ be an oriented hypersurface immersed in the unit standard
sphere. Now we denote by $\vec{x}$ its position vector of point $p$
of $N^{n}$ and by $\vec{n}$ the unit normal vector field of $N^{n}$
in $S^{n+1}(1)$. Then we can define its \lq\lq Gauss map\rq\rq by
\begin{equation*}
{\mathcal G}:N^{n}\ni{p} \longmapsto
[\vec{x}(p)\wedge\vec{n}(p)]\cong [\vec{x}(p)+\sqrt{-1}\vec{n}(p)]
\in\widetilde{\mathrm {Gr}}_{2}({\mathbf R}^{n+2}) \cong
Q_{n}({\mathbf C}).
\end{equation*}
Here $[\vec{x}(p)\wedge \vec{n}(p)]$ denotes an oriented $2$-plane
in ${\mathbf R}^{n+2}$ spanned by two vectors $\vec{x}(p)$ and
$\vec{n}(p)$. Then the fundamental fact is that ${\mathcal G}$ is a
Lagrangian immersion (\cite{Palmer94},\cite{Palmer97}). Moreover we
can observe that

\begin{proposition}\label{HypsurfDeform}
Let $F:N^{n}\rightarrow S^{n+1}(1)$ be a smooth immersion of an
$n$-dimensional oriented smooth manifold $N^{n}$ into the
$(n+1)$-dimensional unit sphere $S^{n+1}(1)$ and
$\mathcal{G}:N^{n}\rightarrow \widetilde{\mathrm{Gr}}_{2}({\mathbf
R}^{n+2})\cong{Q_{n}({\mathbf C})}$ be the Gauss map of $F$.
\begin{enumerate}
\renewcommand{\labelenumi}{(\arabic{enumi})}
\item
If $F_{t}:N^{n}\rightarrow S^{n+1}(1){\ }(\vert{t}\vert<{c})$ is a
smooth family of smooth immersions with $F=F_{0}$, then a smooth
family of the Gauss maps $\mathcal{G}_{t}:N^{n}\rightarrow
\widetilde{\mathrm{Gr}}_{2}({\mathbf R}^{n+2})\cong{Q_{n}({\mathbf
C})}$ of $F_{t}{\ }(\vert{t}\vert<{c})$ is a Hamiltonian deformation
of $\mathcal{G}$.
\item
Suppose that $\varphi_{t}:N^{n}\rightarrow
\widetilde{\mathrm{Gr}}_{2}({\mathbf R}^{n+2})\cong{Q_{n}({\mathbf
C})}{\ } (\vert{t}\vert<{c})$ is a Hamiltonian deformation of the
Gauss map $\varphi_{0}=\mathcal{G}$. If $N^{n}$ is compact or
$\{\varphi_{t}\}$ is compactly supported, then there exists a
positive real number $\delta<c$ and a smooth family of smooth
immersions $F_{t}:N^{n}\rightarrow S^{n+1}(1){\
}(\vert{t}\vert<\delta)$ with $F_{0}=F$ such that the Gauss map of
$F_{t}$ coincides with $\varphi_{t}$ for each $t$ with
$\vert{t}\vert<\delta$.
\end{enumerate}
\end{proposition}

\begin{proof}
(1) For each point $p\in{N}$, we denote by $\vec{x}_{t}(p)$ the
position vector of the point $F(p)$ for the immersion
$F_{t}:N^{n}\rightarrow S^{n+1}(1){\ }(\vert{t}\vert<c)$. Let
$\vec{n}_{t}$ denote the unit normal vector field of the immersion
$F_{t}$ compatible with the orientations. Then
$\{\vec{x}_{t},\vec{n}_{t}\}$ defines an orthonormal frame field of
the pull-back vector bundle $\mathcal{G}_{t}^{\ast}\mathcal{V}$
defined on the whole $N^{n}$, which is parallel with respect to the
induced flat connection. Thus for all $t$ the induced flat
connections have trivial holonomy, particularly same holonomy.
Therefore by Theorem
\ref{HamiltonianDeformationsIsomonodromyDeformation} the family
$\mathcal{G}_{t}:N^{n}\rightarrow Q_{n}({\mathbf C}){\
}(\vert{t}\vert<\delta)$ is a Hamiltonian deformation with
${\mathcal G}_{0}={\mathcal G}$.

(2) Assume that $\varphi_{t}:N^{n}\rightarrow Q_{n}({\mathbf C})$ is
a Hamiltonian deformation of $\varphi_{0}={\mathcal G}$. Since the
induced flat $U(1)$-connection in $\varphi_{0}^{-1}{\mathcal
L}={\mathcal G}^{-1}{\mathcal L}$ has trivial holonomy, by Theorem
\ref{HamiltonianDeformationsIsomonodromyDeformation} each induced
flat $U(1)$-connection in $\varphi_{t}^{-1}{\mathcal L}$ also has
trivial holonomy. Thus for each $t$ we can choose smoothly a
parallel orthonormal frame field $\{v_{1}^{t},v_{2}^{t}\}$ in
$\varphi_{t}^{\ast}\mathcal{V}$ on $N$ with
$\{v_{1}^{0},v_{2}^{0}\}=\{{\mathbf x},{\mathbf n}\}$ at $t=0$. Then
a smooth family of smooth maps
$F_{t}:N^{n}\ni{p}\longmapsto{v_{1}^{t}(p)}\in{S^{n+1}(1)}$
satisfies $F_{0}=F$ and hence there is a sufficiently small
$\delta>0$ such that each $F_{t}{\ }(\vert{t}\vert<\delta)$ is an
immersion into $S^{n+1}(1)$ whose Gauss map coincides with
$\varphi_{t}$.
\end{proof}
A local weaker version of Proposition \ref{HypsurfDeform} was stated
also in \cite{Palmer97}. B.~Palmer (\cite{Palmer97}) showed the mean
curvature form formula as follows :

\begin{equation}
\alpha_{H} = d\left(\mathrm{Im} \left(
\log{\prod^{n}_{i=1}(1+\sqrt{-1}\kappa_{i})} \right)\right),
\end{equation}
where $\kappa_{i}{\ }(i=1,\cdots,n)$ denotes the principal
curvatures of $N^{n}\subset S^{n+1}(1)$.

In case $n=2$, since $(1+\sqrt{-1}\kappa_{1})(1+\sqrt{-1}\kappa_{2})
=1-K_{N}+\sqrt{-1}H_{N}$, we see that for any minimal surface
$N^{2}\subset S^{3}(1)$, its Gauss map ${\mathcal G}:N^{2}
\longrightarrow \widetilde{\mathrm{Gr}}_{2}({\mathbf R}^{4}) \cong
Q_{2}({\mathbf C})\cong S^{2}\times S^{2}$ is a minimal Lagrangian
immersion. This case has been investigated by many authors. In the
next section we shall discuss the case of general $n$ when all
principal curvatures $\kappa_{i}$ are constant.

For a point $[V]$ of $\widetilde{\mathrm{Gr}}_{2}({\mathbf
R}^{n+2})$, the fiber of the vector bundle $\mathcal{V}$ at $[V]$ is
given by
$$
(\mathcal{V})_{[V]}=\{([V],v){\ }\vert{\ }v\in{V}\}
$$
where $[V]$ is an oriented $2$-dimensional vector subspace of
${\mathbf R}^{n+2}$ considered as a point of
$\widetilde{\mathrm{Gr}}_{2}({\mathbf R}^{n+2})$. For each
$[V]\in\widetilde{\mathrm{Gr}}_{2}({\mathbf R}^{n+2})$, let $V\oplus
V^{\perp}={\mathbf R}^{n+2}$ be the decomposition into $V$ and its
orthogonally complementary subspace $V^{\perp}$. We have the
identification:
\begin{equation*}
T_{[V]}\widetilde{\mathrm{Gr}}_{2}({\mathbf R}^{n+2})\cong
{\mathrm{Hom}}(V,V^{\perp}).
\end{equation*}
The standard complex structure ${\mathcal J}$ of
$\widetilde{\mathrm{Gr}}_{2}({\mathbf R}^{n+2})$ is defined by
\begin{equation*}
[{\mathcal J}(T)](v):=T(jv)
\end{equation*}
for each $T\in{\mathrm{Hom}}(V,V^{\perp})\cong
T_{[V]}\widetilde{\mathrm{Gr}}_{2}({\mathbf R}^{n+2})$ and each
$v\in{V}$, where $j$ is a rotation of $\pi/2$ on the oriented
$2$-dimensional vector space $V$.

Let $\varphi : L\rightarrow \widetilde{\mathrm{Gr}}_{2}({\mathbf
R}^{n+2})$ be a Lagrangian immersion of an $n$-dimensional connected
smooth manifold $L$ and let $g_{L}$ denote a Riemannian metric on
$L$ induced by $\varphi$ from the standard Riemannian metric of
$\widetilde{\mathrm{Gr}}_{2}({\mathbf R}^{n+2})$.

Let $\pi_{\varphi}:\varphi^{-1}\mathcal{V}\longrightarrow L$ be the
pull-back vector bundle over $L$ with induced flat $U(1)$-connection
$\varphi^{-1}\nabla$. For each
$(p_{0},v_{0})\in{\varphi^{-1}\mathcal{V}}$, where
$([V_{0}],v_{0})\in{\mathcal{V}}_{\varphi(p_{0})}$ with
$\varphi(p_{0})=[V_{0}]\in\widetilde{\mathrm{Gr}}_{2}({\mathbf
R}^{n+2})$ and a unit vector $v_{0}\in{V_{0}}\subset{\mathbf
R}^{n+2}$, there is a unique maximal connected integral manifold
$\tilde{N}$ through $(p_{0},v_{0})$ in $\varphi^{-1}\mathcal{V}$ of
the horizontal distribution with respect to the flat connection
$\varphi^{-1}\nabla$. Then $\pi_{\varphi}:\tilde{N}\rightarrow{L}$
is a smooth covering map with the deck transformation group
$\rho(\pi_{1}(L))\subset{U(1)}$, where
$\rho:\pi_{1}(L)\rightarrow{U(1)}$ is the holonomy homomorphism of
the flat connection $\varphi^{-1}\nabla$. We define a smooth map $F$
by
\begin{equation}
F:\tilde{N}\ni{(x,v)}\longmapsto v\in{S^{n+1}(1)}.
\end{equation}
We should note that the map $F$ is not always an immersion and thus
the image $F(\tilde{N})$ does not necessarily give a hypersurface in
$S^{n+1}(1)$.

We express $\varphi$ as
\begin{equation*}
\varphi:L\ni{p}\longmapsto \varphi(p)=[V_{p}]\in
{\widetilde{\mathrm{Gr}}_{2}({\mathbf R}^{n+2})},
\end{equation*}
where $V_{p}$ is an oriented $2$-dimensional vector subspace of
${\mathbf R}^{n+2}$ representing the point $\varphi(p)$. The
differential of $\varphi$
\begin{equation*}
(d\varphi)_{p}:T_{p}L\longrightarrow
T_{[V_{p}]}\widetilde{\mathrm{Gr}}_{2}({\mathbf R}^{n+2})\cong
{\mathrm{Hom}}(V_{p},V_{p}^{\perp})
\end{equation*}
is an injective isometric linear map. For each vector $v\in{V_{p}}$,
we introduce a linear operator
\begin{equation*}
B_{v}:T_{p}L\rightarrow{V_{p}}^{\perp}\subset{\mathbf R}^{n+2}
\end{equation*}
defined by
\begin{equation*}
B_{v}(X):=[(d\varphi)_{p}(X)](v)
\end{equation*}
for each $X\in{T_{p}L}$.

\begin{proposition}
The linear operators $B_{v}{\ }(v\in{V_{p}})$ have the following
properties :
\begin{equation}\label{LagrCond0}
{\ }^{t}B_{v}\circ{B_{jv}} = {\ }^{t}B_{jv}\circ{B_{v}} \quad\mbox{
for each }v\in{V_{p}},
\end{equation}
and
\begin{equation}\label{IsometCond0}
{\ }^{t}B_{v}\circ{B_{v}} + {\ }^{t}B_{jv}\circ{B_{jv}}=\mathrm{Id}
\quad\mbox{ for each }v\in{V_{p}}.
\end{equation}
\end{proposition}

\begin{proof}
It follows from the Lagrangian condition of $\varphi$ that
\begin{equation*}\label{LagrCond2}
\begin{split}
&\langle{[(d\varphi)_{p}(X)](jv_{1}),[(d\varphi)_{p}(Y)](v_{1})}\rangle
+
\langle{[(d\varphi)_{p}(X)](jv_{2}),[(d\varphi)_{p}(Y)](v_{2})}\rangle\\
=&
\langle{[(d\varphi)_{p}(X)](jv_{1}),[(d\varphi)_{p}(Y)](v_{1})}\rangle
-
\langle{[(d\varphi)_{p}(X)](v_{1}),[(d\varphi)_{p}(Y)](jv_{1})}\rangle\\
=&0
\end{split}
\end{equation*}
for each $X,Y\in T_{p}L$, where $\{v_{1},v_{2}=jv_{1}\}$ is an
orthonormal basis of $V_{p}$ compatible with its orientation.
Thus
\begin{equation}\label{LagrCond3}
\langle{B_{jv_{1}}(X),B_{v_{1}}(Y)}\rangle =
\langle{B_{v_{1}}(X),B_{jv_{1}}(Y)}\rangle
\end{equation}
for each $X,Y\in T_{p}L$. The condition (\ref{LagrCond3}) is equal
to
\begin{equation*}\label{LagrCond4}
{\ }^{t}B_{v_{1}}\circ{B_{jv_{1}}}={\
}^{t}B_{jv_{1}}\circ{B_{v_{1}}}.
\end{equation*}
On the other hand, the isometric condition of $\varphi$ implies
\begin{equation*}\label{IsometCond2}
\begin{split}
&\langle{[(d\varphi)_{p}(X)](v_{1}),[(d\varphi)_{p}(Y)](v_{1})}\rangle
+
\langle{[(d\varphi)_{p}(X)](jv_{1}),[(d\varphi)_{p}(Y)](jv_{1})}\rangle\\
=& \langle{X,Y}\rangle
\end{split}
\end{equation*}
for each $X,Y\in T_{p}L$.
So we have
\begin{equation*}\label{IsometCond3}
{\ }^{t}B_{v_{1}}\circ{B_{v_{1}}} + {\
}^{t}B_{jv_{1}}\circ{B_{jv_{1}}} =\mathrm{Id}.
\end{equation*}
\end{proof}

Moreover,
by using (\ref{LagrCond0}) and (\ref{IsometCond0}) we get

\begin{lemma}
For any
$v=\cos{\theta}\cdot{v_{1}}+\sin{\theta}\cdot{jv_{1}}\in{V_{p}}$,
$$
{\ }^{t}B_{v}\circ B_{v} =\cos{2\theta}{\ }(^{t}B_{v_{1}}\circ
B_{v_{1}}) +\sin{2\theta}{\ }(^{t}B_{jv_{1}}\circ B_{v_{1}})
+\sin^{2}\theta{\ }\mathrm{I}.
$$
\end{lemma}
If we suppose that $0\not=X\in\mathrm{Ker}B_{v_{1}}\subset{T_{p}L}$,
then we have
\begin{eqnarray*}
(^{t}B_{v}\circ B_{v})X &=& \cos{2\theta}{\ }(^{t}B_{v_{1}}\circ
B_{v_{1}})X +\sin{2\theta}{\ }(^{t}B_{jv_{1}}\circ B_{v_{1}})X
+\sin^{2}\theta{\ }X\\
&=&\sin^{2}\theta{\ }X.
\end{eqnarray*}
Thus we see that if $0<\theta<\pi$, then $X\notin\mathrm{Ker}B_{v}$.
Note that
$\Vert{B_{v}(X)}\Vert=\vert{\sin\theta}\vert\cdot\Vert{X}\Vert$.
More strongly we obtain

\begin{lemma}
$\mathrm{Ker}B_{v}\perp \mathrm{Ker}B_{v_{1}}$ for each
$\theta\mbox{ with }0<\theta<\pi$.
\end{lemma}

\begin{proof}
For each $X\in\mathrm{Ker}B_{v}$ and each
$Y\in\mathrm{Ker}B_{v_{1}}$, we have
\begin{equation*}
\begin{split}
0=&\langle(^{t}B_{v}\circ B_{v})X,Y\rangle\\
=&\cos{2\theta}{\ }\langle(^{t}B_{v_{1}}\circ B_{v_{1}})X,Y\rangle
+\sin{2\theta}{\ }\langle(^{t}B_{jv_{1}}\circ B_{v_{1}})X,Y\rangle
+\sin^{2}\theta{\ }\langle{X,Y}\rangle\\
=&\cos{2\theta}{\ }\langle(^{t}B_{v_{1}}\circ B_{v_{1}})X,Y\rangle
+\sin{2\theta}{\ }\langle(^{t}B_{v_{1}}\circ B_{jv_{1}})X,Y\rangle
+\sin^{2}\theta{\ }\langle{X,Y}\rangle\\
=&\cos{2\theta}{\ }\langle{B_{v_{1}}X,B_{v_{1}}Y}\rangle
+\sin{2\theta}{\ }\langle{B_{jv_{1}}X,B_{v_{1}}Y}\rangle
+\sin^{2}\theta{\ }\langle{X,Y}\rangle\\
=&\sin^{2}\theta{\ }\langle{X,Y}\rangle.
\end{split}
\end{equation*}
\end{proof}

Therefore we obtain

\begin{lemma}\label{InjBv}
Except for finitely many $v\in{V_{p}}\cap{S^{n+1}(1)}$,
$B_{v}:T_{p}L\rightarrow{V_{p}}^{\perp}$ is a linear isomorphism.
\end{lemma}

For such a vector $v$, there is a unique maximal connected integral
manifold $\tilde{N}$ through $(p,v)$ in $\varphi^{-1}\mathcal{V}$ of
the horizontal distribution with respect to the flat connection
$\varphi^{-1}\nabla$. The map $F:\tilde{N}\ni{(x,w)}\longmapsto
w\in{S^{n+1}(1)}$ is an immersion in a neighborhood of $(p,v)$ and
$F(\tilde{N})$ gives an immersed oriented hypersurface $S^{n+1}(1)$
around $(p,v)$.

\section{Lagrangian submanifolds in complex hyperquadrics obtained as
Gauss images of isoparametric hypersurfaces in spheres}\label{sec:3}

Now suppose that $N^{n}$ is a compact oriented hypersurface in
$S^{n+1}(1)$ with constant principal curvatures, the so called
\lq\lq isoparametric hypersurface\rq\rq. By M\"unzner's result
(\cite{Muenzner1},\cite{Muenzner2}), $N^{n}$ is real algebraic in
the sense that it is defined by a certain homogeneous real algebraic
equation (Cartan-M\"unzner polynomial) and the number $g$ of
distinct principal curvatures must be $g=1,2,3,4 \mbox{ or } 6$.
Then the \lq\lq{Gauss image}\rq\rq of a minimal Lagrangian immersion
${\mathcal G}:N^{n}\rightarrow Q_{n}({\bold C})$ is a compact
minimal Lagrangian submanifold $L={\mathcal G}(N^{n})=N^{n}/{\bold
Z}_{g}$ embedded in ${Q_{n}({\bold C})}$ obtained as the quotient
space of $N^{n}$ by a free action of a finite cyclic group ${\bold
Z}_{g}$ of order $g$. We remark that $g=1\mbox{ or }2$ if and only
if ${\mathcal G}:N^{n}\rightarrow Q_{n}({\bold C})$ is a totally
geodesic Lagrangian immersion. All isoparametric hypersurfaces in
spheres are classified into homogeneous ones, which are given as
principal orbits of compact group actions on spheres with
cohomogeneity $1$, and non-homogneous ones, which were discovered
first by H.~Ozeki-M.~Takeuchi (\cite{Ozeki-TakeuchiI},
\cite{Ozeki-TakeuchiII}), and developed by
D.~Ferus-H.~Karcher-H.~F.~M\"unzner (\cite{FKM}) and recently
T.~Cecil-Q.-S.~Chi-G.~R.~Jensen (\cite{Cecil-Chi-Jensen}). Concerned
with the homogeneity, we can observe

\begin{proposition}\label{homogeneity}
An isoparametric hypersurface $N^{n}\subset S^{n+1}(1)$ is {\it
homogeneous}, i.e., an orbit of a compact connected Lie subgroup
$K\subset{SO(n+2)}$, if and only if its Gauss image ${\mathcal
G}(N^{n})$ is a {\it homogeneous} Lagrangian submanifold in
$Q_{n}({\bold C})$.
\end{proposition}

The part of \lq\lq only if \rq\rq is trivial. Here we give a proof
for the part of \lq\lq if \rq\rq. Assume that $\mathcal{G}(N^{n})$
is homogeneous, that is, a Lagrangian orbit
$\mathcal{G}(N^{n})=K\cdot[V_{0}]$ through a point $[V_{0}]$ of a
compact connected Lie subgroup $K$ of $SO(n+2)$. Then $g_{L}$ is a
$K$-invariant Riemannian metric on $L$. In order to prove that $N$
is also homogeneous, we analyze the $K$-equivariance of the bundle
homomorphism
$$
B:\mathcal{V}\otimes TL \ni v\otimes X \longmapsto B_{v}(X)\in
\mathcal{V}^{\perp}.
$$
For $p\in{N}$, let
$[V_{\mathcal{G}(p)}]=\mathcal{G}(p)$ denote an oriented
$2$-dimensional vector subspace of ${\bold R}^{n+2}$ spanned by
${\bold x}(p)$ and ${\bold n}(p)$ with the orientation determined by
$\{{\bold x}(p),{\bold n}(p)\}$.

Let $\{{\bold e}_{1},\cdots,{\bold e}_{n}\}$ be a local orthonormal
frame field compatible with the orientation of $N$ around $p$ such
that $h({\bold e}_{i},{\bold e}_{j})=\kappa_{i}\delta_{ij}$
$(i,j=1,\cdots,n)$, or equivalently $A({\bold
e}_{i})=\kappa_{i}{\bold e}_{i}$ $(i=1,\cdots,n)$, where $h$ and $A$
denote the second fundamental form and the shape operator or the
Weingarten map of the hypersurface $N^{n}$ in $S^{n+1}(1)$,
respectively. We denote by $\{k_{1},\cdots,k_{g}\}$ the distinct
principal curvatures of $N^{n}$ in $S^{n+1}(1)$.

For each $X\in{T_{p}N}$ and each $v=\cos\phi{\ }{\bold
x}(p)+\sin\phi{\ }{\bold n}(p)\in{V_{\mathcal{G}(p)}}$,
$$
B_{v}((d\mathcal{G})_{p}(X))= [(d\mathcal{G})_{p}(X)](v)= \cos\phi
X-\sin\phi A(X).
$$
Since
$$
B_{v}((d\mathcal{G})_{p}({\bold e}_{i}))
=(\cos{\phi}-{\kappa_{i}}\sin{\phi}){\bold e}_{i}
\quad(i=1,\cdots,n),
$$
the linear endomorphism
$^{t}B_{v}\circ{B}_{v}:T_{\mathcal{G}(p)}L\rightarrow{T_{\mathcal{G}(p)}L}$
satisfies
$$ (^{t}B_{v}\circ{B}_{v})((d\mathcal{G})_{p}({\bold
e}_{i}))
=\frac{(\cos{\phi}-{\kappa_{i}}\sin{\phi})^{2}}{1+\kappa_{i}^{2}}
(d\mathcal{G})_{p}({\bold e}_{i})\quad(i=1,\cdots,n)
$$
Note that
$$
\frac{(\cos{\phi}-{\kappa_{i}}\sin{\phi})^{2}}{1+\kappa_{i}^{2}} =
\frac{(\cos{\phi}-{\kappa_{j}}\sin{\phi})^{2}}{1+\kappa_{j}^{2}}
$$
for all $\phi\in{\bold R}$ if and only if $\kappa_{i}=\kappa_{j}$.
The eigenspace decomposition of $A$ at $p$
$$
T_{p}N=(E_{1})_{p}\oplus\cdots\oplus(E_{g})_{p}
$$
corresponds the eigenspace decomposition of the symmetric linear
endomorphism $^{t}B_{v}\circ{B}_{v}$ at $\mathcal{G}(p)$
$$
T_{\mathcal{G}(p)}L=(E^{\prime}_{1})_{\mathcal{G}(p)}
\oplus\cdots\oplus(E^{\prime}_{g})_{\mathcal{G}(p)}
$$
via the linear isomorphism $(d\mathcal{G})_{p}:T_{p}N\longrightarrow
T_{\mathcal{G}(p)}L$. Then on each
$(E^{\prime}_{i})_{\mathcal{G}(p)}$ we have
$$
^{t}B_{v}\circ{B_{v}}=
\frac{(\cos\phi-k_{i}\sin\phi)^{2}}{1+k_{i}^{2}}
\mathrm{Id}_{(E^{\prime}_{i})_{_{\mathcal{G}(p)}}}.
$$

\begin{lemma}\label{K-equiv}
Let $p_{0},p\in{N^{n}}$ and $a\in{K}\subset{SO(n+2)}$. Assume that
$\mathcal{G}(p)=a\mathcal{G}(p_{0})$. Then we have
\begin{equation*}
({\bold x}(p),{\bold n}(p)) = \pm(a{\bold x}(p_{0}),a{\bold
n}(p_{0}))
\begin{pmatrix}
\cos(m\pi/g)&-\sin(m\pi/g)\\
\sin(m\pi/g)&\cos(m\pi/g)
\end{pmatrix}.
\end{equation*}
for some integer $m$.
\end{lemma}

\begin{proof}
For arbitrary $v_{0}\in{V_{\mathcal{G}(p_0)}}$, we set
$v=av_{0}\in{V_{\mathcal{G}(p)}}\subset{\bold R^{n+2}}$. We express
$v_{0}$ as $v_{0}=(\cos{\theta}){\bold
x}(p_{0})+(\sin{\theta}){\bold n}(p_{0})$, and we also express $v$
as $v=(\cos\phi){\bold x}(p)+(\sin\phi){\bold n}(p)$. Since
$v=av_{0}= (\cos{\theta})a{\bold x}(p_{0})+(\sin{\theta})a{\bold
n}(p_{0})$, setting
\begin{equation*}
({\bold x}(p),{\bold n}(p)) = (a{\bold x}(p_{0}),a{\bold n}(p_{0}))
\begin{pmatrix}
\cos\psi&-\sin\psi\\
\sin\psi&\cos\psi
\end{pmatrix}
\end{equation*}
for some $\psi\in{\bold R}$, we get
$\phi+\psi\equiv{\theta}\mod{2\pi}$. By the $K$-equivariance we have
\begin{equation}\label{K-equivariance}
(^{t}B_{v}\circ{B_{v}})\circ(da)_{\mathcal{G}(p_{0})}
=(da)_{\mathcal{G}(p_{0})}\circ(^{t}B_{v_{0}}\circ{B_{v_{0}}}).
\end{equation}
Through a linear isomorphism $(da)_{\mathcal{G}(p_{0})}:
T_{\mathcal{G}(p_{0})}L \rightarrow T_{\mathcal{G}(p_{0})}L$, the
direct sum decomposition
$$
T_{\mathcal{G}(p_{0})}L =(E^{\prime}_{1})_{\mathcal{G}(p_{0})}
\oplus\cdots\oplus (E^{\prime}_{g})_{\mathcal{G}(p_{0})}
$$
is mapped to a direct sum decomposition
$$
T_{\mathcal{G}(p)}L=
(da)_{\mathcal{G}(p_{0})}(E^{\prime}_{1})_{\mathcal{G}(p_{0})}
\oplus\cdots\oplus
(da)_{\mathcal{G}(p_{0})}(E^{\prime}_{g})_{\mathcal{G}(p_{0})}.
$$
Thus by \eqref{K-equivariance} on each
$(da)_{\mathcal{G}(p_{0})}(E^{\prime}_{j})_{\mathcal{G}(p_{0})}$ we
have
$$
^{t}B_{v}\circ{B_{v}}=
\frac{(\cos\theta-k_{j}\sin\theta)^{2}}{1+k_{j}^{2}}
\mathrm{Id}_{(da)_{\mathcal{G}(p_{0})}}(E^{\prime}_{j})_{\mathcal{G}(p_{0})}.
$$
Hence for each $i(=1,2,\cdots,g)$, there is uniquely
$j(=1,2,\cdots,g)$ such that
\begin{equation*}
(E^{\prime}_{i})_{\mathcal{G}(p)} =
(da)_{\mathcal{G}(p_{0})}(E^{\prime}_{j})_{\mathcal{G}(p_{0})}.
\end{equation*}
Since we have
\begin{equation*}
\begin{split}
\frac{(\cos\theta-k_{j}\sin\theta)^{2}} {1+k_{j}^{2}}
&=\frac{(\cos\phi-k_{i}\sin\phi)^{2}}{1+k_{i}^{2}}\\
&=
\frac{(\cos(\theta-\psi)-k_{i}\sin(\theta-\psi))^{2}}{1+k_{i}^{2}}
\end{split}
\end{equation*}
for all $\theta\in{\bold R}$,
it implies the equation
$(1+k_{i}k_{j})\sin{\psi}=(k_{i}-k_{j})\cos{\psi}$ for each $i$.
Suppose that $\sin{\psi}\not=0$. Then $k_{i}-k_{j}\not=0$ and thus
it becomes
\begin{equation*}
\cot\psi=\frac{1+k_{i}k_{j}}{k_{i}-k_{j}}=\cot(\beta_{i}-\beta_{j})
\qquad\text{ for all }i.
\end{equation*}
Here set $k_{i}=\cot{\beta_{i}}{\ }(0<{\beta_{i}}<\pi)$. Since we
know from the theory of isoparametric hypersurfaces
(\cite{Muenzner1}) that $\beta_{i}-\beta_{j}\in(\pi/g){\bold Z}$, we
obtain that $\psi\in(\pi/g){\bold Z}+\pi{\bold Z}$.
\end{proof}

By the continuation argument on $N^{n}$ and $K$, it follows from
Lemma \ref{K-equiv} that $N^{n}=K\cdot{p_{0}}$ and thus $N^{n}$ is
homogeneous. We complete the proof of Proposition \ref{homogeneity}.

\bigskip

By virtue of the results of Hsiang-Lawson
(\cite{Hsiang-Lawson1971}), Takagi-Takahashi
(\cite{Takagi-Takahashi1972}), we know that all homogeneous
isoparametric hypersurfaces $N^{n}\subset S^{n+1}(1)$ can be
obtained as principal orbits of compact Riemannian symmetric pairs
$(U,K)$ of rank $2$ :
\smallskip
\begin{center}
Table $2$
\end{center}
\begin{center}
\begin{tabular}
{|c|c|c|c|c|c|} \hline $g$ &Type &$(U,K)$ &$\dim{N}$ &$m_{1},m_{2}$
&$N=K/K_{0}$
\\
\hline $1$ &$S^{1}\times$ &$(S^{1}\times SO(n+2),SO(n+1))$ &${n}$
&${n}$ &$S^{n}$
\\
${\ }$ &$\mathrm{BDII}$ &$(n\geq{1})$ &${\ }$ &${\ }$ &${\ }$
\\
\hline $2$ &$\mathrm{BDII}$ &$(SO(p+2)\times SO(n+2-p),$ &$n$
&$p,n-p$ &$S^{p}\times S^{n-p}$
\\
${\ }$ &$\times\mathrm{BDII}$ &$SO(p+1)\times SO(n+1-p))$ &${\ }$
&${\ }$ &${\ }$
\\
${\ }$ &{\ } &$(1\leq{p}\leq{n-1})$ &${\ }$ &${\ }$ &${\ }$
\\
\hline $3$ &$\mathrm{AI}_{2}$ &$(SU(3),SO(3))$ &$3$ &$1,1$
&$\frac{SO(3)}{{\bold Z}_{2}+{\bold Z}_{2}}$
\\
\hline $3$ &${\frak a}_{2}$ &$(SU(3)\times SU(3),SU(3))$ &$6$ &$2,2$
&$\frac{SU(3)}{T^{2}}$
\\
\hline $3$ &$\mathrm{AII}_{2}$ &$(SU(6),Sp(3))$ &$12$ &$4,4$
&$\frac{Sp(3)}{Sp(1)^{3}}$
\\
\hline $3$ &$\mathrm{EIV}$ &$(E_{6},F_{4})$ &$24$ &$8,8$
&$\frac{F_{4}}{Spin(8)}$
\\
\hline $4$ &${\frak b}_{2}$ &$(SO(5)\times SO(5),SO(5))$ &$8$ &$2,2$
&$\frac{SO(5)}{T^{2}}$
\\
\hline $4$ &$\mathrm{AIII}_{2}$ &$(SU(m+2),S(U(m)\times U(2)))$
&$4m-2$ &$2,$ &$\frac{S(U(m)\times U(2))}{SU(m-2)\times{T^{2}}}$
\\
${\ }$ &{\ } &$(m\geq{2})$ &${\ }$ &$2m-3$ &${\ }$
\\
\hline $4$ &$\mathrm{BDI}_{2}$ &$(SO(m+2),SO(m)\times SO(2))$
&$2m-2$ &$1,$ &$\frac{SO(m)\times SO(2)}{SO(m-2)\times{\bold
Z}_{2}}$
\\
${\ }$ &{\ } &$(m\geq{3})$ &${\ }$ &$m-2$ &${\ }$
\\
\hline $4$ &$\mathrm{CII}_{2}$ &$(Sp(m+2),Sp(m)\times Sp(2))$
&$8m-2$ &$4,$ &$\frac{Sp(m)\times Sp(2)}{Sp(m-2)\times{Sp(1)^{2}}}$
\\
${\ }$ &{\ } &$(m\geq{2})$ &${\ }$ &$4m-5$ &${\ }$
\\
\hline $4$ &$\mathrm{DIII}_{2}$ &$(SO(10),U(5))$ &$18$ &$4,5$
&$\frac{U(5)}{{SU(2)}\times{SU(2)}\times{T^{1}}}$
\\
\hline $4$ &$\mathrm{EIII}$ &$(E_{6},Spin(10)\cdot{T})$ &$30$ &$6,9$
&$\frac{Spin(10)\cdot{T}}{SU(4)\cdot{T}}$
\\
\hline $6$ &${\frak g}_{2}$ &$(G_{2}\times G_{2},G_{2})$ &$12$
&$2,2$ &$\frac{G_{2}}{T^{2}}$
\\
\hline $6$ &$\mathrm{G}$ &$(G_{2},SO(4))$ &$6$ &$1,1$
&$\frac{SO(4)}{{\bold Z}_{2}+{\bold Z}_{2}}$
\\
\hline
\end{tabular}
\end{center}
\medskip
Here $m_{1}$, $m_{2}$ denote the multiplicities of the principal
curvatures of $N^{n}$.
Let ${\frak u}={\frak k}+{\frak p}$ be the canonical decomposition
of ${\frak u}$ as a symmetric Lie algebra of a symmetric pair
$(U,K)$ and ${\frak a}$ be a maximal abelian subspace of ${\frak
p}$. Let $B_{\frak u}({\ },{\ })$ denote the Killing-Cartan form of
the Lie algebra ${\frak u}$ if $(U,K)$ is not of type
$S^{1}\times\mathrm{BDII}$, or the direct sum of the $(-1)$-times
standard inner product of $\sqrt{-1}{\bold R}$ and the
Killing-Cartan form of ${\frak u}$ if $(U,K)$ is of type
$S^{1}\times\mathrm{BDII}$. Define the $\mathrm{Ad}(U)$-invariant
inner product $\langle{\ },{\ }\rangle_{\frak u}$ of ${\frak u}$ by
$\langle{X},{Y}\rangle_{\frak u}:=-B_{\frak u}(X,Y)$ for each
$X,Y\in{\frak u}$. The vector space ${\frak p}$ is identified with
the Euclidean space ${\bold R}^{n+2}$ with respect to the inner
product $\langle{\ },{\ }\rangle_{\frak u}$. The $(n+1)$-dimensional
unit sphere $S^{n+1}(1)$ in ${\frak p}$ is defined as
\begin{equation*}
S^{n+1}(1):=\{X\in{\frak p}{\ }\vert{\ } \Vert{X}\Vert_{\frak
u}^{2}=\langle{X},{X}\rangle_{\frak u}=1\}.
\end{equation*}
The isotropy linear action $\mathrm{Ad}_{\frak p}$ of $K$ on ${\frak
p}$ and thus $S^{n+1}(1)$ induces the group action of $K$ on
$\widetilde{\mathrm{Gr}}_{2}({\frak p})\cong Q_{n}({\bold C})$. For
each {\it regular} element $H$ of ${\frak a}\cap{S^{n+1}(1)}$, we
get a homogeneous isoparametric hypersurface in the unit sphere
\begin{equation*}
N^{n}=(\mathrm{Ad}_{\frak p}K)H\subset S^{n+1}(1)\subset{\bold
R}^{n+2} \cong{\frak p}.
\end{equation*}
Its Gauss image is
\begin{equation*}
{\mathcal G}(N^{n}) =K\cdot[{\frak a}] =[(\mathrm{Ad}_{\frak
p}K){\frak a}]\subset \widetilde{\mathrm{Gr}}_{2}({\frak p})\cong
Q_{n}({\bold C}).
\end{equation*}
Here $N$ and ${\mathcal G}(N^{n})$ have homogeneous space
expressions $N\cong K/K_{0}$ and ${\mathcal G}(N^{n})\cong
K/K_{[{\frak a}]}$, where we set
\begin{equation*}
\begin{split}
&K_{0}:=\{k\in{K}{\ }\vert{\ }\mathrm{Ad}_{\frak p}(k)(H)=H\}\\
&\quad{\ }=\{k\in{K}{\ }\vert{\ }\mathrm{Ad}_{\frak p}(k)(H)=H
\text{ for each }H\in{\frak a}\},\\
&K_{\frak a}:=\{k\in{K}{\ }\vert{\ } \mathrm{Ad}_{\frak p}(k)({\frak
a})
={\frak a}\},\\
&K_{[{\frak a}]}:=\{k\in{K_{\frak a}}{\ }\vert{\ }
\mathrm{Ad}_{\frak p}(k):{\frak a}\longrightarrow{\frak a} \text{
preserves the orientation of }{\frak a}\}.
\end{split}
\end{equation*}
The deck transformation group of the covering map ${\mathcal G}:
N\rightarrow {\mathcal G}(N^{n})$ is equal to $K_{[{\frak
a}]}/K_{0}=W(U,K)/{\bold Z}_{2}\cong{\bold Z}_{g}$, where
$W(U,K)=K_{\frak a}/K_{0}$ is the Weyl group of $(U,K)$.

The moment map $\mu$ of the action of $K$ on $Q_{n}({\bold C})$
induced by the adjoint action of $K$ on ${\frak p}$ is given as
follows :
\begin{equation}\label{MomentMap(U,K)}
\mu:Q_{n}({\bold C})\cong\widetilde{\mathrm{Gr}}_{2}({\frak p}) \ni
[{\bold a}+\sqrt{-1}{\bold b}]=[V]\longmapsto -[{\bold a},{\bold
b}]\in{\frak k}\cong{\frak k}^{\ast}
\end{equation}
where $\{{\bold a},{\bold b}\}$ is an orthonormal basis of
$V\subset{\frak p}$ compatible with its orientation. Then we obtain
\begin{equation}
{\mathcal G}(N^{n})=\mu^{-1}(0).
\end{equation}
For each $[V]\in{\widetilde{\mathrm{Gr}}_{2}({\frak p})}$, the
square norm of the moment map
\begin{equation}
\Vert{\mu([V])}\Vert^{2}=\Vert{[{\bold a},{\bold b}]}\Vert^{2}
\end{equation}
is equal to the sectional curvature of compact symmetric space
$(U/K,g_{U})$ for a $2$-plane $V$ (cf.~\cite{Helgason78},
\cite{KNI-II}). Here $g_{U}$ denotes the invariant Riemannian metric
on $U/K$ induced by the $\mathrm{Ad}(U)$-invariant inner product
$\langle{\ },{\ }\rangle_{\frak u}$ of ${\frak u}$.

Let $\Sigma(U,K)$ be the set of (restricted) roots of $({\frak
u},{\frak k})$ and $\Sigma^{+}(U,K)$ be its subset of positive
roots. Note that each $\gamma\in{\Sigma(U,K)}$ is an ${\bold
R}$-linear function $\gamma:{\frak a}\rightarrow{\sqrt{-1}{\bold
R}}$.
We have the root decomposition of ${\frak k}$ and ${\frak p}$ as
follows (cf.~\cite{Takagi-Takahashi1972}):
\begin{equation}\label{RestRootDecomp}
{\frak k} ={\frak k}_{0}+\sum_{\gamma\in{\Sigma^{+}(U,K)}} {\frak
k}_{\gamma},\quad {\frak p} ={\frak
a}+\sum_{\gamma\in{\Sigma^{+}(U,K)}}{\frak p}_{\gamma},
\end{equation}
where
\begin{equation*}
\begin{split}
{\frak k}_{0}:
=&\{X\in{\frak k}{\ }\vert{\ }[X,{\frak a}]\subset{\frak a}\}\\
=&\{X\in{\frak k}{\ }\vert{\ }[X,H]=0\quad
\text{ for each }H\in{\frak a}\},\\
{\frak k}_{\gamma} :=&\{X\in{\frak k}{\ }\vert{\
}(\mathrm{ad}H)^{2}X=(\gamma(H))^{2}X
\text{ for each }H\in{\frak a}\},\\
{\frak p}_{\gamma} :=&\{Y\in{\frak p} {\ }\vert{\ }
(\mathrm{ad}H)^{2}Y=(\gamma(H))^{2}Y \text{ for each }H\in{\frak
a}\}.
\end{split}
\end{equation*}
For each $\gamma\in{\Sigma^{+}(U,K)}$, set $m(\gamma):=\dim{{\frak
k}_{\gamma}} =\dim{{\frak p}_{\gamma}}$. Define
\begin{equation}\label{M&Aperp}
{\frak m}:= \sum_{\gamma\in{\Sigma^{+}(U,K)}}{\frak k}_{\gamma}\quad
\text{ and }\quad {\frak a}^{\perp}:=
\sum_{\gamma\in{\Sigma^{+}(U,K)}}{\frak p}_{\gamma}.
\end{equation}
We can choose orthonormal bases of ${\frak m}$ and ${\frak
a}^{\perp}$ with respect to $\langle{\ },{\ }\rangle_{\frak u}$
\begin{equation}\label{StdBasis1}
\{X_{\gamma,i}{\ }\vert{\
}\gamma\in{\Sigma^{+}(U,K)},i=1,2,\cdots,m(\gamma)\}
\end{equation}
and
\begin{equation}\label{StdBasis2}
\{Y_{\gamma,i}{\ }\vert{\
}\gamma\in{\Sigma^{+}(U,K)},i=1,2,\cdots,m(\gamma)\}
\end{equation}
such that
\begin{equation}\label{StdBasis}
[H,X_{\gamma,i}]=\sqrt{-1}\gamma(H)Y_{\gamma,i},\quad
[H,Y_{\gamma,i}]=-\sqrt{-1}\gamma(H)X_{\gamma,i}
\end{equation}
for each $H\in{\frak a}$.

The following condition \eqref{CenterCond} is used in the
classification of Lagrangian orbits in complex hyperquadrics in the
next section.

\begin{lemma}\label{CenterCondition}
Assume that $(U,K)$ is a compact Riemannian symmetric pair in
$\mathrm{Table}{\ }2$. Then $(U,K)$ satisfies the condition
\begin{equation}\label{CenterCond}
\{0\}\not={\frak c}({\frak k})\subset{\frak m}
\end{equation}
if and only if $(U,K)$ is one of the following pairs :
\begin{enumerate}
\renewcommand{\labelenumi}{(\alph{enumi})}
\item
$(S^{1}\times{SO(3)},SO(2)){\ } [S^{1}\times\mathrm{BDII},n=1]$.
\item
$(SO(3)\times SO(3),SO(2)\times SO(2)) {\
}[\mathrm{BDII}\times\mathrm{BDII},n=2]$.
\item
$(SO(3)\times SO(n+1),SO(2)\times SO(n)){\ }(n\geq{3}) {\
}[\mathrm{BDII}\times\mathrm{BDII},p=1\text{ or }n-1]$.
\item
$(SO(m+2),SO(2)\times SO(m)){\ }(n=2m-2,m\geq{3}) {\
}[\mathrm{BDI}_{2}]$.
\end{enumerate}
\end{lemma}

We shall prove Lemma \ref{CenterCondition}. According to Table $2$,
if $K$ has the nontrivial center, then $(U,K)$ must be $(S^{1}\times
SO(3),SO(2))$, $(SO(3)\times SO(n+1),SO(2)\times SO(n)){\
}(n\geq{2})$ or an irreducible Hermitian symmetric pair $(U,K)$ of
compact type and of rank $2$. It is obvious that $(U,K)=(S^{1}\times
SO(3),SO(2))$ and $(U,K)=(SO(3)\times SO(n+1),SO(2)\times SO(n)){\
}(n\geq{2})$ satisfy the condition \eqref{CenterCond}. Now we assume
that $(U,K)$ is an irreducible Hermitian symmetric pair of compact
type and of rank $2$. Note that $\dim{\frak c}({\frak k})=1$. We
refer \cite{Helgason78} for results from the theory of Hermitian
symmetric Lie algebras. Let ${\frak u}={\frak k}+{\frak p}$ be the
canonical decomposition of ${\frak u}$ with respect to a Hermitian
symmetric pair $(U,K)$ and ${\frak a}$ be a maximal abelian subspace
of ${\frak p}$. Let ${\frak g}={\frak k}+\sqrt{-1}{\frak p}$ denote
its noncompact dual in the complexification ${\frak u}^{\bold
C}={\frak g}^{\bold C}$ and ${\frak g}^{\bold C}={\frak k}^{\bold
C}+{\frak p}^{\bold C}$ their complexification. Let ${\frak
c}({\frak k})$ be the center of ${\frak k}$ and ${\frak h}$ be a
maximal abelian subalgebra of ${\frak k}$ and ${\frak u}$ so that
${\frak c}({\frak k})\subset{\frak h}\subset{\frak k} \subset{\frak
u}$.
Let $\Delta$ denote the set of all roots of ${\frak g}^{\bold C}$
with respect to ${\frak h}^{\bold C}$ and $\Delta^{+}$ denote the
set of all positive roots in $\Delta$ relative to a suitable linear
order. $\alpha\in\Delta$ is called {\it compact} (resp.~{\it
noncompact}) if $\alpha=0$ (resp.~$\alpha\not=0$) on ${\frak
c}({\frak k})$. We take the root decompositions with respect to
${\frak h}^{\bold C}$ :
\begin{equation}
{\frak u}^{\bold C}={\frak h}^{\bold C}
+\sum_{\alpha\in\Delta}{\frak g}^{\alpha},\quad {\frak k}^{\bold
C}={\frak h}^{\bold C} +\sum_{\alpha\text{:cpt}}{\frak g}^{\alpha},
\quad {\frak p}^{\bold C}= \sum_{\alpha\text{:noncpt}}{\frak
g}^{\alpha}.
\end{equation}
Set $\Delta_{\frak c}:=\{\alpha\in\Delta{\ }\vert{\ }
\alpha\not=0\text{ on }{\frak c}({\frak k})\}$ and
$Q_{+}:=\Delta^{+}\cap\Delta_{\frak c}$.
Then there exist strongly orthogonal roots
$\gamma_{1},\gamma_{2}\in{Q_{+}}$ such that
\begin{equation}
{\frak a}^{\bold C}= \sum^{2}_{i=1}{\bold
C}(X_{\gamma_{i}}+X_{-\gamma_{i}})
\end{equation}
and
\begin{equation}\label{a}
{\frak a}= \sum^{2}_{i=1}{\bold
R}\sqrt{-1}(X_{\gamma_{i}}+X_{-\gamma_{i}}).
\end{equation}
Here we use $X_{\alpha}\in{\frak g}^{\alpha}$ such that for each
$\alpha\in\Delta$,
\begin{equation}\label{StdBasisWRTDelta}
\begin{split}
&X_{\alpha}-X_{-\alpha},\sqrt{-1}(X_{\alpha}+X_{-\alpha})\in{\frak u},\\
&[X_{\alpha},X_{-\alpha}]=\frac{2}{\alpha(H_{\alpha})}H_{\alpha}.
\end{split}
\end{equation}
Define an abelian subalgebra ${\frak h}_{\frak a}^{\bold C}$ of
${\frak k}_{0}^{\bold C}$ by
\begin{equation}
{\frak h}_{\frak a}^{\bold C} :=\{H\in{\frak h}^{\bold C}{\ }\vert{\
} \gamma_{i}(H)=\langle{H_{\gamma_{i}}},{H}\rangle=0 {\ }(i=1,2)\}.
\end{equation}

\begin{lemma}\label{ZinCenter}
Let $Z\in{\frak c}({\frak k})$ be a nonzero element. Then
$Z\in{\frak m}$ if and only if $Z\perp{{\frak h}_{\frak a}}$, if and
only if $Z\in\sum^{2}_{i=1}{\bold R}H_{\gamma_{i}}$.
\end{lemma}

\begin{proof}
Suppose that $X_{\alpha}\in{\frak g}^{\alpha}$ with compact $\alpha$
satisfies
\begin{equation*}
[X_{\alpha},X_{\gamma_{i}}+X_{-\gamma_{i}}] \in{\frak a}^{\bold
C}\quad(i=1,2).
\end{equation*}
Then the condition
\begin{equation*}
0= [X_{\alpha},X_{\gamma_{i}}+X_{-\gamma_{i}}]=
[X_{\alpha},X_{\gamma_{i}}] +[X_{\alpha},X_{-\gamma_{i}}] \in{\ }
{\frak g}^{\alpha+\gamma_{i}} +{\frak g}^{\alpha-\gamma_{i}}
\end{equation*}
implies that $[X_{\alpha},X_{\gamma_{i}}]=0$ and
$[X_{\alpha},X_{-\gamma_{i}}]=0$. We see that $\alpha+\gamma_{i}$
and $\alpha-\gamma_{i}$ are not roots. Thus we obtain the expression
\begin{equation}\label{K0C}
{\frak k}_{0}^{\bold C}={\frak h}_{\frak a}^{\bold C}
+\sum_{\alpha}{\frak g}^{\alpha},
\end{equation}
where $\alpha\in\Delta$ runs all compact roots such that
$\alpha\pm\gamma_{i}$ are not roots for all $i=1,2$. Hence from
\eqref{K0C} we have
\begin{equation}\label{MC}
{\frak m}^{\bold C}= \left(\sum^{2}_{i=1}{\bold
C}H_{\gamma_{i}}\right) +\sum_{\alpha}{\frak g}^{\alpha},
\end{equation}
where $\alpha\in\Delta$ runs all compact roots such that
$\alpha+\gamma_{i}$ or $\alpha-\gamma_{i}$ is a root for some
$i=1,2$. And we have an orthogonal direct sum decomposition
\begin{equation}\label{H}
{\frak h}= \left(\sum^{2}_{i=1}{\bold R}H_{\gamma_{i}}\right)
\oplus{\frak h}_{\frak a}.
\end{equation}
From \eqref{MC} and \eqref{H} we obtain Lemma \ref{ZinCenter}.
\end{proof}


By using Lemma \ref{ZinCenter}, we determine Hermitian symmetric
pairs $(U,K)$ of compact type and of rank $2$ satisfying the
condition \eqref{CenterCond} as follows. We use the table of root
systems for complex simple Lie algebras in \cite{Bourbaki}.

Let $\{\alpha_{1},\cdots,\alpha_{l}\}$ be the fundamental root
system of $\Delta$. We express the highest root $\tilde{\alpha}$ of
$\Delta$ as
\begin{equation*}
\tilde{\alpha}=m_{1}\alpha_{1}+\cdots+m_{l}\alpha_{l}.
\end{equation*}
Define a basis $\{\xi_{1},\cdots,\xi_{l}\}$ of ${\frak h}$ by
\begin{equation*}
\alpha_{i}(\xi_{j})=\sqrt{-1}\langle{H_{\alpha_{i}}},{\xi_{j}}\rangle
=\sqrt{-1}\delta_{ij}\quad(i,j=1,\cdots,l).
\end{equation*}
The fundamental weight system $\{\Lambda_{1},\cdots,\Lambda_{l}\}$
is defined by
\begin{equation*}
2{\ }\frac{\langle{\Lambda_{j}},{\alpha_{i}}\rangle}
{\langle{\alpha_{i}},{\alpha_{i}}\rangle} =-2\sqrt{-1}{\
}\frac{{\alpha_{i}}(H_{\Lambda_{j}})}
{\langle{\alpha_{i}},{\alpha_{i}}\rangle}
=\delta_{ij}\quad(i,j=1,\cdots,l).
\end{equation*}
Note that
$$
\xi_{j}=\frac{2H_{\Lambda_{j}}}{\langle{\alpha_{j}},{\alpha_{j}}\rangle}
\quad(j=1,2,\cdots,\l).
$$
In our case there is a subset $I=\{i_{0}\}$ of $\{1,2,\cdots,l\}$
with $m_{i_{0}}=1$ such that
\begin{equation}
{\frak k}^{\bold C}={\frak h}^{\bold C}
+\sum_{n_{i_{0}}(\alpha)=0}{\frak g}^{\alpha}, \quad {\frak
p}^{\bold C}= \sum_{n_{i_{0}}(\alpha)\not=0}{\frak g}^{\alpha},
\end{equation}
where we define $n_{i}(\alpha)$ as
$\alpha=\sum^{l}_{i=1}n_{i}(\alpha)\alpha_{i}\in\Delta$ (all $n_{i}$
are nonnegative integers or all $n_{i}$ are nonpositive integers).
Then we have
\begin{equation*}
{\frak c}({\frak k})={\bold R}\xi_{i_{0}} \subset{\bold
C}\xi_{i_{0}} ={\frak c}({\frak k}^{\bold C}) ={\frak c}({\frak
k})^{\bold C}.
\end{equation*}
Since $\Delta_{\frak{c}}=\{\alpha\in\Delta{\ }\vert{\ }
n_{i_{0}}(\alpha)=1\}$, note that
$Q_{+}=\Delta^{+}\cap\{\alpha\in\Delta{\ }\vert{\ } n_{i_{0}}=1\}$.

$\mathrm{AIII}_{2}$ : The case when $(U,K)=(SU(m+2),S(U(m)\times
U(2))){\ }(m\geq{2})$ is described as
\begin{equation*}
\begin{split}
&\tilde{\alpha}=\alpha_{1}+\alpha_{2}+\cdots+\alpha_{m+1}
=\sqrt{-1}(\varepsilon_{1}-\varepsilon_{m+2}),{\ }
H_{\tilde{\alpha}}=e_{1}-e_{m+2},\\
&i_{0}=2,{\ }
\alpha_{2}=\sqrt{-1}(\varepsilon_{2}-\varepsilon_{3}),{\ }
H_{\alpha_{2}}=e_{2}-e_{3},\\
&\xi_{2} =(e_{1}+e_{2})-\frac{2}{m+2}\sum^{m+2}_{j=1}e_{j}.
\end{split}
\end{equation*}
The strongly orthogonal roots are given by $\gamma_{1}=\alpha_{2},{\
}\gamma_{2}=\tilde{\alpha}\in{Q_{+}}$. Therefore
\begin{equation*}
\xi_{2}\notin {\bold R}H_{\alpha_{2}}+{\bold R}H_{\tilde{\alpha}} =
{\bold R}(e_{2}-e_{3}) +{\bold R}(e_{1}-e_{m+2})
\end{equation*}
if and only if $m\geq{3}$. If $m=2$, then
\begin{equation*}
\begin{split}
\xi_{2} =&\frac{1}{2}(e_{1}+e_{2}) -\frac{1}{2}(e_{3}+e_{4})
=\frac{1}{2}(e_{2}-e_{3})
+\frac{1}{2}(e_{1}-e_{4})\\
\in& {\bold R}H_{\alpha_{2}}+{\bold R}H_{\tilde{\alpha}} = {\bold
R}(e_{2}-e_{3}) +{\bold R}(e_{1}-e_{4}).
\end{split}
\end{equation*}
Note that $\mathrm{AIII}_{2}{\ }(m=2)\cong\mathrm{BDI}_{2}{\
}(m=4)$.

$\mathrm{BII}_{2}$ : The case when $(U,K)=(SO(m+2),SO(2)\times
SO(m)), m=2l-1, l\geq{2}$ is described as
\begin{equation*}
\begin{split}
&\tilde{\alpha}=\alpha_{1}+2\alpha_{2}+\cdots+2\alpha_{l}
=\sqrt{-1}(\varepsilon_{1}+\varepsilon_{2}),{\ }
H_{\tilde{\alpha}}=e_{1}+e_{2}.\\
&i_{0}=1,{\ }
\alpha_{1}=\sqrt{-1}(\varepsilon_{1}-\varepsilon_{2}),{\ }
H_{\alpha_{1}}=e_{1}-e_{2},\\
&\xi_{1}=e_{1}.
\end{split}
\end{equation*}
The strongly orthogonal roots are given by $\gamma_{1}=\alpha_{1},{\
} \gamma_{2}=\tilde{\alpha}\in{Q_{+}}$. Thus
\begin{equation*}
\xi_{1}=e_{1}\in {\bold R}H_{\alpha_{1}}+{\bold
R}H_{\tilde{\alpha}}= {\bold R}e_{1}+{\bold R}e_{2}.
\end{equation*}
Hence the condition \eqref{CenterCond} is satisfied in this case.

$\mathrm{DII}_{2}$ : The case when $(U,K)=(SO(m+2),SO(2)\times
SO(m)), m=2l, l\geq{2}$ is described as
\begin{equation*}
\begin{split}
&\tilde{\alpha}=\alpha_{1}+2\alpha_{2}+\cdots+2\alpha_{l-2}
+\alpha_{l-1}+\alpha_{l}
=\sqrt{-1}(\varepsilon_{1}+\varepsilon_{2}),{\ }
H_{\tilde{\alpha}}=e_{1}+e_{2}.\\
&i_{0}=1,{\ }
\alpha_{1}=\sqrt{-1}(\varepsilon_{1}-\varepsilon_{2}),{\ }
H_{\alpha_{1}}=e_{1}-e_{2},\\
&\xi_{1}=e_{1}.
\end{split}
\end{equation*}
The strongly orthogonal roots are given by $\gamma_{1}=\alpha_{1},{\
} \gamma_{2}=\tilde{\alpha}\in{Q_{+}}$.

\begin{equation*}
\xi_{1}=e_{1}\in {\bold R}H_{\alpha_{1}}+{\bold
R}H_{\tilde{\alpha}}= {\bold R}e_{1}+{\bold R}e_{2}.
\end{equation*}
Therefore the condition \eqref{CenterCond} is satisfied in this
case.

$\mathrm{DIII}_{2}$ : The case when $(U,K)=(SO(10),U(5))$ is
described as
\begin{equation*}
\begin{split}
&\tilde{\alpha}=\alpha_{1}+2\alpha_{2}+2\alpha_{3}
+\alpha_{4}+\alpha_{5}
=\sqrt{-1}(\varepsilon_{1}+\varepsilon_{2}),{\ }
H_{\tilde{\alpha}}=e_{1}+e_{2}.\\
&i_{0}=4,{\ }
\alpha_{4}=\sqrt{-1}(\varepsilon_{4}-\varepsilon_{5}),{\ }
H_{\alpha_{4}}=e_{4}-e_{5},\\
&\xi_{4}=\frac{1}{2} (e_{1}+e_{2}+e_{3}+e_{4}-e_{5}).
\end{split}
\end{equation*}
The strongly orthogonal roots are given by $\gamma_{1}=\alpha_{4},
\gamma_{2}=\tilde{\alpha}\in{Q_{+}}$. Hence
\begin{equation*}
\xi_{4}\notin {\bold R}H_{\alpha_{4}}+{\bold R}H_{\tilde{\alpha}}.
\end{equation*}
Therefore the condition \eqref{CenterCond} is not satisfied in this
case.

$\mathrm{EIII}$ : The case when $(U,K)=(E_{6},Spin(10)\cdot{T})$ is
described as
\begin{equation*}
\begin{split}
\tilde{\alpha}&= \alpha_{1}+2\alpha_{2}+2\alpha_{3}
+3\alpha_{4}+2\alpha_{5}+\alpha_{6}\\
&= \frac{1}{2}\sqrt{-1}(\varepsilon_{1}+\varepsilon_{2}
+\varepsilon_{3}+\varepsilon_{4} +\varepsilon_{5}-\varepsilon_{6}
-\varepsilon_{7}+\varepsilon_{8}),\\
H_{\tilde{\alpha}}&= \frac{1}{2}(e_{1}+e_{2}+e_{3}+e_{4}
+e_{5}-e_{6}-e_{7}+e_{8}).\\
i_{0}&=1,{\ } \alpha_{1}=\frac{1}{2}\sqrt{-1}
(\varepsilon_{1}+\varepsilon_{8})
-\frac{1}{2}\sqrt{-1}(\varepsilon_{2}+\varepsilon_{3}
+\varepsilon_{4}+\varepsilon_{5} +\varepsilon_{6}+\varepsilon_{7}),
\\
H_{\alpha_{1}}&= \frac{1}{2}(e_{1}+e_{8})
-\frac{1}{2}(e_{2}+e_{3}+e_{4}+e_{5}+e_{6}+e_{7}),\\
\xi_{1}&=\frac{2}{3} (e_{8}-e_{7}-e_{6}).
\end{split}
\end{equation*}
The strongly orthogonal roots are given by $\gamma_{1}=\alpha_{1},
\gamma_{2}=\tilde{\alpha}\in{Q_{+}}$. Hence
\begin{equation*}
\xi_{1}\notin {\bold R}H_{\alpha_{1}}+{\bold R}H_{\tilde{\alpha}}.
\end{equation*}
Therefore the condition \eqref{CenterCond} is not satisfied in this
case.
We complete the proof of Lemma \ref{CenterCondition}.

Using \eqref{a}, we set
\begin{equation}\label{H1H2}
\begin{split}
2H_{1}:&=\sqrt{-1}(X_{\gamma_{1}}+X_{-\gamma_{1}})
-\sqrt{-1}(X_{\gamma_{2}}+X_{-\gamma_{2}})\in{\frak a},\\
2H_{2}:&=\sqrt{-1}(X_{\gamma_{1}}+X_{-\gamma_{1}})
+\sqrt{-1}(X_{\gamma_{2}}+X_{-\gamma_{2}})\in{\frak a}.
\end{split}
\end{equation}
Then the center ${\frak c}({\frak k})$ of ${\frak k}$ is given as
follows :

\begin{lemma}\label{[JH_{2},H_{2}]}
Suppose that $[\mathrm{BII}_{2}]$ $(U,K)=(SO(m+2),SO(2)\times
SO(m)), m=2l-1,l\geq{2}$ or $[\mathrm{DII}_{2}]$ :
$(U,K)=(SO(m+2),SO(2)\times SO(m)), m=2l,l\geq{2}$. Then we have
\begin{equation}
0\not=[JH_{2},H_{2}]\in{\frak c}(\frak k).
\end{equation}
\end{lemma}

\begin{proof}
Using \eqref{StdBasisWRTDelta} we compute
\begin{equation*}
\begin{split}
4[JH_{2},H_{2}]
=& -\alpha_{1}(\xi_{1})
[X_{\alpha_{1}}-X_{-\alpha_{1}},X_{\alpha_{1}}+X_{-\alpha_{1}}]\\
& -\alpha_{1}(\xi_{1}) [X_{\alpha_{1}}-X_{-\alpha_{1}},
X_{\tilde{\alpha}}+X_{-\tilde{\alpha}}]\\
& -\tilde{\alpha}(\xi_{1}) [X_{\tilde{\alpha}}-X_{-\tilde{\alpha}},
X_{\alpha_{1}}+X_{-\alpha_{1}}]\\
& -\tilde{\alpha}(\xi_{1}) [X_{\tilde{\alpha}}-X_{-\tilde{\alpha}},
X_{\tilde{\alpha}}+X_{-\tilde{\alpha}}]\\
=& -2\alpha_{1}(\xi_{1})
[X_{\alpha_{1}},X_{-\alpha_{1}}]\\
& -\alpha_{1}(\xi_{1}) ([X_{\alpha_{1}},X_{-\tilde{\alpha}}]
-[X_{-\alpha_{1}},X_{\tilde{\alpha}}])
\\
& -\tilde{\alpha}(\xi_{1}) ([X_{\tilde{\alpha}},X_{-\alpha_{1}}]
-[X_{-\tilde{\alpha}},X_{\alpha_{1}}])
\\
& -2\tilde{\alpha}(\xi_{1})
[X_{\tilde{\alpha}},X_{-\tilde{\alpha}}]\\
=& -\alpha_{1}(\xi_{1})\cdot
2\frac{2}{\alpha_{1}(H_{\alpha_{1}})}H_{\alpha_{1}}
-\tilde{\alpha}(\xi_{1})\cdot
2\frac{2}{\tilde{\alpha}(H_{\tilde{\alpha}})}H_{\tilde{\alpha}}\\
=& -2(H_{\alpha_{1}}+H_{\tilde{\alpha}})
=-4\xi_{1}\in{\frak c}({\frak k}).
\end{split}
\end{equation*}
In this computation we use the following properties of this type
(cf.~\cite{Bourbaki}): (i) $\alpha_{1}+\tilde{\alpha}\not\in\Delta$.
(ii) $\alpha_{1}-\tilde{\alpha}\not\in\Delta$,
$-\alpha_{1}+\tilde{\alpha}\not\in\Delta$. (iii)
$\alpha_{1}(H_{\alpha_{1}})
=\sqrt{-1}\langle{H_{\alpha_{1}}},{H_{\alpha_{1}}}\rangle=2\sqrt{-1}$,
$\tilde{\alpha}(H_{\tilde{\alpha}})
=\sqrt{-1}\langle{H_{\tilde{\alpha}}},{H_{\tilde{\alpha}}}\rangle
=2\sqrt{-1}$,
$\alpha_{1}(\xi_{1})=\sqrt{-1}\langle{H_{\alpha_{1}}},{\xi_{1}}\rangle
=\sqrt{-1}$, $\tilde{\alpha}(\xi_{1})
=\sqrt{-1}\langle{H_{\tilde{\alpha}}},{\xi_{1}}\rangle=\sqrt{-1}$.
\end{proof}

\section{Classification of homogeneous Lagrangian submanifolds
in complex hyperquadrics}\label{sec:4}

\subsection{Lagrangian orbits and moment maps
in complex hyperquadrics}

Suppose that $K^{\prime}\cdot[V_{0}]\subset{Q_{n}({\bold C})}$ is a
Lagrangian orbit of a compact connected Lie subgroup $K^{\prime}$ of
$SO(n+2)$ through a point $[V_{0}]\in{Q_{n}({\bold C})}$. Let
$\mu_{K^{\prime}} : Q_{n}({\bold C})\rightarrow{\frak k}^{\prime}$
denote the moment map for the induced Hamiltonian group action of
$K^{\prime}$ on $Q_{n}({\bold C})$. By virtue of Lemma \ref{InjBv}
there is a unit vector $v\in{V_{0}}$ such that the orbit
$K^{\prime}\cdot{v}\subset{S^{n+1}(1)}$ through $v$ under the group
$K^{\prime}$ is a homogeneous isoparametric hypersurface embedded in
$S^{n+1}(1)$. Thus by Hsiang-Lawson's theorem
(\cite{Hsiang-Lawson1971}) there is a Riemannian symmetric pair
$(U,K)$ of compact type and of rank $2$ with compact connected $K$
such that ${\frak p}={\bold R}^{n+2}$,
$K^{\prime}\subset\mathrm{Ad}_{\frak p}(K)$, and
$K^{\prime}\cdot{v}=\mathrm{Ad}_{\frak p}(K)v$, where ${\frak
u}={\frak k}+{\frak p}$ denotes the canonical decomposition of the
symmetric Lie algebra $({\frak u},{\frak k})$ of $(U,K)$. Then we
can show that $K^{\prime}\cdot{[V_{0}]}=K\cdot[V_{0}]$ as follows:

According to the classification due to T.~Asoh (\cite{TAsoh1981}),
if $K^{\prime}\not=\mathrm{Ad}_{\frak p}(K)$, then one of the
following cases happens :
\begin{enumerate}
\renewcommand{\labelenumi}{(\roman{enumi})}
\item
$K^{\prime}$ is semisimple, i.e., ${\frak c}({\frak
k}^{\prime})=\{0\}$.
\item
$K^{\prime}$ is not semisimple and
$K^{\prime}=K^{\prime}_{ss}\times{C(K^{\prime})}$ where
$K^{\prime}_{ss}$ is semisimple such that $K^{\prime}_{ss}$ acts on
$S^{n+1}(1)$ with cohomogeneity $1$.
\item
$K^{\prime}$ is not semisimple and
$K^{\prime}=K^{\prime}_{ss}\times{C(K^{\prime})}$ where
$K^{\prime}_{ss}$ is semisimple such that $K^{\prime}_{ss}$ does NOT
act on $S^{n+1}(1)$ with cohomogeneity $1$. In this case
\begin{enumerate}
\renewcommand{\labelenumi}{(\alph{enumi})}
\item
$K^{\prime}=SO(2)\times G_{2}$, $n=12$,
$(U,K)=(SO(9),SO(2)\times SO(7))$.
\item
$K^{\prime}=SO(2)\times Spin(7)$, $n=14$,
$(U,K)=(SO(10),SO(2)\times SO(8))$.
\end{enumerate}
\end{enumerate}

In case (i), since ${\frak c}({\frak k}^{\prime})=\{0\}$,
$K^{\prime}$ has only one Lagrangian orbit $K^{\prime}\cdot{[V_{0}]}
=\mu_{K^{\prime}}^{-1}(0)=\mathcal{G}(K^{\prime}\cdot{v})
=\mathcal{G}(K\cdot{v})=\mu_{K}^{-1}(0)=K\cdot{[V_{0}]}$ on
$Q_{n}({\bold C})$.

In case (ii), since $K^{\prime}_{ss}\cdot{[V_{0}]}\subset
K^{\prime}\cdot{[V_{0}]}
=(K^{\prime}_{ss}\times{C(K^{\prime})})\cdot{[V_{0}]}$, we see that
$K^{\prime}_{ss}\cdot{[V_{0}]}\subset Q_{n}({\bold C})$ is an
isotropic orbit. Thus
$K^{\prime}_{ss}\cdot{[V_{0}]}\subset\mu_{K^{\prime}_{ss}}^{-1}(0)$.
On the other hand
$\mu_{K^{\prime}_{ss}}^{-1}(0)=\mathcal{G}(K^{\prime}_{ss}\cdot{v})
=\mathcal{G}(K^{\prime}\cdot{v})=\mathcal{G}(K\cdot{v})
=K^{\prime}_{ss}\cdot\mathcal{G}({v})$ is a Lagrangian orbit. Hence
$K^{\prime}_{ss}\cdot{[V_{0}]}=K^{\prime}_{ss}\cdot\mathcal{G}({v})
=\mathcal{G}(K^{\prime}_{ss}\cdot{v})=\mathcal{G}(K^{\prime}\cdot{v})
=\mathcal{G}(K\cdot{v})=K^{\prime}\cdot{[V_{0}]}=K\cdot{[V_{0}]}
\subset Q_{n}({\bold C})$.

In case (iii), we may assume one of the following :
\begin{enumerate}
\renewcommand{\labelenumi}{(\alph{enumi})}
\item
$K^{\prime}=SO(2)\times G_{2}$, $n=12$,
$(U,K)=(SO(9),SO(2)\times SO(7))$ and $K_{0}={\bold Z}_{2}\times SO(5)$,
$K^{\prime}_{0}=K^{\prime}\cap{K_{0}}={\bold Z}_{2}\times SU(2)$.
\item
$K^{\prime}=SO(2)\times Spin(7)$, $n=14$, $(U,K)=(SO(10),SO(2)\times
SO(8))$ and $K_{0}={\bold Z}_{2}\times SO(6)$,
$K^{\prime}_{0}=K^{\prime}\cap{K_{0}}={\bold Z}_{2}\times SU(3)$.
\end{enumerate}
Since $K^{\prime}\cdot{v}=K\cdot{v}$, we can assume that $v\in{\frak
a}$. We take the orthogonal direct sum ${\frak k}={\frak
k}^{\prime}+{{\frak k}^{\prime}}^{\perp}$. And we have ${\frak
k}_{0}=\{\xi\in{\frak k}{\ }\vert{\ }[\xi,v]=0\}$. Set ${\frak
k}^{\prime}_{0}={\frak k}_{0}\cap{\frak k}^{\prime}$ and we take the
orthogonal direct sum ${\frak k}_{0}={\frak k}^{\prime}_{0}+{{\frak
k}^{\prime}_{0}}^{\perp}$. Then we obtain ${{\frak
k}^{\prime}}^{\perp}={{\frak k}^{\prime}_{0}}^{\perp}$. Let
$\{v,jv\}$ be an orthonormal basis of $V_{0}$ compatible with its
orientation.
Since $\mu_{K^{\prime}}([V_{0}])=-[v,jv]_{{\frak k}^{\prime}}
\in{\frak c}({\frak k}^{\prime})$,
we can express
$[v,jv]$ as
\begin{equation*}
[v,jv]=\eta+\zeta,
\end{equation*}
where $\eta\in{\frak c}({\frak k}^{\prime})$ and
$\zeta\in{{\frak k}^{\prime}}^{\perp}$.
Since $v\in{\frak a}$ and
${{\frak k}^{\prime}}^{\perp}=
{{\frak k}^{\prime}_{0}}^{\perp}\subset{\frak k}_{0}$,
we have
\begin{equation*}
[v,[v,jv]]=[v,\eta],
\end{equation*}
and so
\begin{equation*}
\begin{split}
\langle[v,[v,jv]],jv\rangle_{\frak u}&=\langle[v,\eta],jv\rangle_{\frak u}\\
-\langle[v,jv],[v,jv]\rangle_{\frak u}&=-\langle\eta,[v,jv]\rangle_{\frak
u}\\
-\langle{\eta,\eta}\rangle_{\frak
u}-\langle{\zeta,\zeta}\rangle_{\frak u}
&=-\langle\eta,\eta\rangle_{\frak u}\\
\langle{\zeta,\zeta}\rangle_{\frak u}&=0.
\end{split}
\end{equation*}
Thus $\zeta=0$ and $\mu_{K}([V_{0}])=-[v,jv]=\eta\in{\frak c}({\frak
k}^{\prime}) ={\frak c}({\frak k})$, that is,
$K\cdot{[V_{0}]}\subset{Q_{n}({\bold C})}$ is an isotropic orbit.
Hence $K^{\prime}\cdot{[V_{0}]}=K\cdot{[V_{0}]}\subset{Q_{n}({\bold
C})}$ is a Lagrangian orbit.

Therefore we may assume that $K=K^{\prime}$ and thus ${\frak
k}={\frak k}^{\prime}$. We choose an orthonormal basis
$\{v_{1},v_{2}\}$ of $V_{0}$ compatible with the orientation of
$[V_{0}]$. By the adjoint action of $K$ on ${\frak p}$ we can assume
that $v_{1}\in{\frak a}\cap{S^{n+1}(1)}$. Following the
decomposition \eqref{RestRootDecomp}, we decompose $v_{2}$ as
$$
v_{2}=v_{2,0}+\sum_{\gamma\in\Sigma^{+}(U,K)}v_{2,\gamma},
$$
where $v_{2,0}\in{\frak a}$ and $v_{2,\gamma}\in{\frak p}_{\gamma}$.
Then it follows from $[v_{1},v_{2,\gamma}]\in{\frak k}_{\gamma}$
that
\begin{equation*}
\begin{split}
{\frak c}({\frak k})\ni\mu([V_{0}])
&=-[v_{1},v_{2}]\\
&=-[v_{1},v_{2,0}+\sum_{\gamma\in\Sigma^{+}(U,K)}v_{2,\gamma}]\\
&=-\sum_{\gamma\in{\Sigma^{+}(U,K)}}[v_{1},v_{2,\gamma}]\in{\frak
m}.
\end{split}
\end{equation*}
Set $\xi_{0}:=\mu([V_{0}])\in{\frak c}({\frak k})\cap{\frak m}$.

In the case when ${\frak c}({\frak k})\cap{\frak m}=\{0\}$, we have
$\mu([V_{0}])=0$ and thus $K\cdot{[V_{0}]}\subset\mu^{-1}(0)$. Hence
since $\mu^{-1}(0)$ is connected by Proposition \ref{KirwanResult}
(1), we obtain $K\cdot{[V_{0}]}=\mu^{-1}(0)=K\cdot{[{\frak a}]}$.

Suppose that ${\frak c}({\frak k})\cap{\frak m}\not=\{0\}$, or
equivalently $\{0\}\not={\frak c}({\frak k})\subset{\frak m}$. By
Lemma \ref{CenterCondition}, $(U,K)$ must be one of the following :
\begin{enumerate}
\renewcommand{\labelenumi}{(\alph{enumi})}
\item $(S^{1}\times SO(3),SO(2))$,
\item $(SO(3)\times SO(3),SO(2)\times SO(2))$,
\item $(SO(3)\times SO(n+1),SO(2)\times SO(n)){\ }(n\geq{3})$,
\item $(SO(m+2),SO(2)\times SO(m)){\ }(n=2m-2,m\geq{3})$.
\end{enumerate}

\begin{lemma}\label{NormMomentMap}
Set
\begin{equation}
s:=\mathrm{Max}\{\Vert{\mu([V])}\Vert^{2}{\ }\vert{\ }
[V]\in\widetilde{\mathrm{Gr}}_{2}({\frak p})\text{ with }
\mu([V])\in{\frak c}({\frak k})\}.
\end{equation}
\begin{enumerate}
\renewcommand{\labelenumi}{(\alph{enumi})}
\item
If $(U,K)=(S^{1}\times SO(3),SO(2))$, then the value $s$ is equal to
the maximum of sectional curvatures of
$(U/K,g_{U})=S^{1}\times(SO(3)/SO(2))$. Moreover
$\Vert{\mu([V])}\Vert^{2}$ attains the maximum if and only if $[V]$
is an oriented $2$-dimensional subspace of ${\frak p}$ tangent to
$SO(3)/SO(2)$.
\item
If $(U,K)=(SO(3)\times SO(3),SO(2)\times SO(2))$, then the value $s$
is equal to the maximum of sectional curvatures of
$(U/K,g_{U})=(SO(3)/SO(2))\times(SO(3)/SO(2))$. Moreover
$\Vert{\mu([V])}\Vert^{2}$ attains the maximum if and only if $[V]$
is an oriented $2$-dimensional subspace of ${\frak p}$ tangent to
$SO(3)/SO(2)$ with larger curvature.
\item
If $(U,K)=(SO(3)\times SO(n+1),SO(2)\times SO(n)){\ }(n\geq{3})$,
then the value $s$ is equal to the sectional curvature of the first
factor $SO(3)/SO(2)$ of
$(U/K,g_{U})=(SO(3)/SO(2))\times(SO(n+1)/SO(n))$. Moreover
$\Vert{\mu([V])}\Vert^{2}$ is equal to the sectional curvature of
the first factor $SO(3)/SO(2)$ if and only if $[V]$ is an oriented
$2$-dimensional subspace of ${\frak p}$ tangent to $SO(3)/SO(2)$.
\item
If $(U,K)=(SO(m+2),SO(2)\times SO(m)){\ }(n=2m-2,m\geq{3})$, then
the value $s$ is equal to the minimum of holomorphic sectional
curvatures of $(U/K,g_{U},J)=SO(m+2)/(SO(2)\times
SO(m))=Q_{m}({\bold C})$. Moreover $\Vert{\mu([V])}\Vert^{2}$ is
equal to the minimum of holomorphic sectional curvatures of
$(U/K,g_{U})$ if and only if $[V]$ is an oriented $2$-dimensional
vector subspace of ${\frak p}$ spanned by a vector $H_{\gamma}$ dual
to a short root $\gamma\in\Sigma(U,K)$ and $JH_{\gamma}$, up to the
adjoint action of $K$ on ${\frak p}$.
\end{enumerate}
\end{lemma}

\begin{proof}

(a) Let $(U,K)=(S^{1}\times SO(3),SO(2))$. Then
\begin{equation*}
\begin{split}
{\frak u}&=\sqrt{-1}{\bold R}\oplus\frak{so}(3),\quad {\frak
k}=\frak{so}(2)={\frak c}({\frak k}) =\Bigl\{
\begin{pmatrix}
0&0&0\\
0&0&c\\
0&-c&0
\end{pmatrix}
{\ }\vert{\ }c\in{\bold R}
\Bigr\}\\
{\frak p}&=\sqrt{-1}{\bold R}\oplus \left\{
\begin{pmatrix}
0&x&y\\
-x&0&0\\
-y&0&0
\end{pmatrix}
{\ }\vert{\ }x,y\in{\bold R}
\right\},\\
{\frak a}&=\sqrt{-1}{\bold R}\oplus \left\{
\begin{pmatrix}
0&x&0\\
-x&0&0\\
0&0&0
\end{pmatrix}
{\ }\vert{\ }x\in{\bold R} \right\}.
\end{split}
\end{equation*}
Set
\begin{equation*}
H_{1}:=(\sqrt{-1},0)\in\sqrt{-1}{{\bold R}\oplus\{0\}}\subset{\frak
a},\quad H_{2}:=
\begin{pmatrix}
0&1&0\\
-1&0&0\\
0&0&0
\end{pmatrix}
\in{\frak a}.
\end{equation*}
In this case we define $J:=\mathrm{ad}_{\frak p}Z$ by
\begin{equation*}
Z=
\begin{pmatrix}
0&0&0\\
0&0&1\\
0&-1&0
\end{pmatrix}.
\end{equation*}
Then we obtain
\begin{equation*}
[\cos\theta H_{1}+\sin\theta JH_{2},H_{2}]
=\sin\theta[JH_{2},H_{2}]=(\sin\theta)Z
\end{equation*}
and
\begin{equation*}
\Vert[\cos\theta H_{1}+\sin\theta JH_{2},H_{2}]\Vert^{2}
=\vert\sin\theta\vert^{2}\Vert{Z}\Vert^{2}.
\end{equation*}
This implies the statement (a).

(b) Let $(U,K)=(SO(3)\times SO(3),SO(2)\times SO(2))$. Then
\begin{equation*}
\begin{split}
{\frak u}&=\frak{so}(3)\oplus\frak{so}(3),\quad {\frak k}=
\frak{so}(2)\oplus\frak{so}(2)={\frak c}({\frak k}),\\
{\frak p}&= \left\{ \left(
\begin{pmatrix}
0&x&y\\
-x&0&0\\
-y&0&0
\end{pmatrix},
\begin{pmatrix}
0&x^{\prime}&y^{\prime}\\
-x^{\prime}&0&0\\
-y^{\prime}&0&0
\end{pmatrix}
\right) {\ }\vert{\ }x,y,x^{\prime},y^{\prime}\in{\bold R}
\right\},\\
{\frak a}&= \left\{ \left(
\begin{pmatrix}
0&x&0\\
-x&0&0\\
0&0&0
\end{pmatrix},
\begin{pmatrix}
0&x^{\prime}&0\\
-x^{\prime}&0&0\\
0&0&0
\end{pmatrix}
\right) {\ }\vert{\ }x,x^{\prime}\in{\bold R} \right\}.
\end{split}
\end{equation*}
Let
\begin{equation*}
v_{1}= \left(
\begin{pmatrix}
0&x_{1}&0\\
-x_{1}&0&0\\
0&0&0
\end{pmatrix},
\begin{pmatrix}
0&x^{\prime}_{1}&0\\
-x^{\prime}_{1}&0&0\\
0&0&0
\end{pmatrix}
\right) \in{\frak a}
\end{equation*}
with $(x_{1})^{2}+(x^{\prime}_{1})^{2}=1$. We may assume that
$x_{1}\geq{0}$ and $x^{\prime}_{1}\geq{0}$. Let
\begin{equation*}
v_{2}= \left(
\begin{pmatrix}
0&x_{2}&y_{2}\\
-x_{2}&0&0\\
-y_{2}&0&0
\end{pmatrix},
\begin{pmatrix}
0&x^{\prime}_{2}&y^{\prime}_{2}\\
-x^{\prime}_{2}&0&0\\
-y^{\prime}_{2}&0&0
\end{pmatrix}
\right)\in{\frak p}
\end{equation*}
with ${x_{2}}^{2}+{y_{2}}^{2}
+{x^{\prime}_{2}}^{2}+{y^{\prime}_{2}}^{2}=1$. Then
\begin{equation*}
{\frak k}={\frak c}({\frak k})\ni [v_{1},v_{2}]= \left( x_{1}y_{2}
\begin{pmatrix}
0&0&0\\
0&0&-1\\
0&1&0
\end{pmatrix},
x^{\prime}_{1}y^{\prime}_{2}
\begin{pmatrix}
0&0&0\\
0&0&-1\\
0&1&0
\end{pmatrix}
\right)\in{\frak c}({\frak k})={\frak k}
\end{equation*}
and thus
\begin{equation*}
\begin{split}
\Vert[{v_{1}},{v_{2}}]\Vert^{2} =&(x_{1}y_{2})^{2}
+(x^{\prime}_{1}y^{\prime}_{2})^{2}\\
=&(x_{1})^{2} +(x^{\prime}_{1})^{2}\leq{1}.
\end{split}
\end{equation*}
Here $\Vert[{v_{1}},{v_{2}}]\Vert^{2}=1$ if and only if
\begin{equation*}
(y_{2})^{2}=1 \text{ or } (y^{\prime}_{2})^{2}=1, \text{ i.e., }{\ }
y_{2}=\pm{1} \text{ or } y^{\prime}_{2}=\pm{1},
\end{equation*}
if and only if
\begin{equation*}
v_{1}=\left(
\begin{pmatrix}
0&1&0\\
-1&0&0\\
0&0&0
\end{pmatrix},
0\right),\quad v_{2}=\left(
\begin{pmatrix}
0&0&\pm{1}\\
0&0&0\\
\mp{1}&0&0
\end{pmatrix},
0\right)
\end{equation*}
\text{ or }
\begin{equation*}
v_{1}=\left(0,
\begin{pmatrix}
0&1&0\\
-1&0&0\\
0&0&0
\end{pmatrix}
\right),\quad v_{2}=\left(0,
\begin{pmatrix}
0&0&\pm{1}\\
0&0&0\\
\mp{1}&0&0
\end{pmatrix}
\right).
\end{equation*}
This implies the statement (b).

(c) Let $(U,K)=(SO(3)\times SO(n+1),SO(2)\times SO(n)){\
}(n\geq{3})$.
\begin{equation*}
\begin{split}
{\frak u}&=\frak{so}(3)\oplus\frak{so}(n+1),{\ } {\frak
k}=\frak{so}(2)\oplus\frak{so}(n),{\ }
{\frak c}({\frak k})=\frak{so}(2)\oplus\{0\},\\
{\frak p}&= \left\{ \left(
\begin{pmatrix}
0&x&y\\
-x&0&0\\
-y&0&0
\end{pmatrix},
\begin{pmatrix}
0&X\\
-^{t}X&0
\end{pmatrix}
\right) {\ }\vert{\ }x,y\in{\bold R},X\in{M(1,n;{\bold R})}
\right\},\\
{\frak a}&= \left\{ \left(
\begin{pmatrix}
0&x&0\\
-x&0&0\\
0&0&0
\end{pmatrix},
\begin{pmatrix}
0&y&0&\cdots&0\\
-y&0&0&\cdots&0\\
0&0&0&\cdots&0\\
\vdots&\vdots&\vdots&\cdots&\vdots\\
0&0&0&\cdots&0
\end{pmatrix}
\right) {\ }\vert{\ }x,y\in{\bold R} \right\}.
\end{split}
\end{equation*}
In this case we define $J:=\mathrm{ad}_{\frak p}Z$ by
\begin{equation*}
Z= \left(
\begin{pmatrix}
0&0&0\\
0&0&-1\\
0&1&0
\end{pmatrix},
0 \right).
\end{equation*}
If we set
\begin{equation*}
H_{1}=\left( 0,
\begin{pmatrix}
0&1&0&\cdots&0\\
-1&0&0&\cdots&0\\
0&0&0&\cdots&0\\
\vdots&\vdots&\vdots&\cdots&\vdots\\
0&0&0&\cdots&0
\end{pmatrix}
\right), \quad H_{2}= \left(
\begin{pmatrix}
0&1&0\\
-1&0&0\\
0&0&0
\end{pmatrix},
0 \right),
\end{equation*}
then
\begin{equation*}
JH_{2}=[Z,H_{2}]= \left(
\begin{pmatrix}
0&0&1\\
0&0&0\\
-1&0&0
\end{pmatrix},
0 \right).
\end{equation*}
Let
\begin{equation*}
{\frak a}\ni v_{1}= \left(
\begin{pmatrix}
0&x&0\\
-x&0&0\\
0&0&0
\end{pmatrix},
\begin{pmatrix}
0&y&0&\cdots&0\\
-y&0&0&\cdots&0\\
0&0&0&\cdots&0\\
\vdots&\vdots&\vdots&\cdots&\vdots\\
0&0&0&\cdots&0
\end{pmatrix}
\right)
\end{equation*}
with $x^{2}+y^{2}=1$. We may assume that $x\geq{0}$. Let
\begin{equation*}
{\frak p}\ni v_{2}= \left(
\begin{pmatrix}
0&x^{\prime}&y^{\prime}\\
-x^{\prime}&0&0\\
-y^{\prime}&0&0
\end{pmatrix},
\begin{pmatrix}
0&X\\
-^{t}X&0
\end{pmatrix}
\right)
\end{equation*}
with ${x^{\prime}}^{2}+{y^{\prime}}^{2}+X{\ }^{t}X=1$. Then
\begin{equation*}
{\frak c}({\frak k})\ni [v_{1},v_{2}] = \left( xy^{\prime}
\begin{pmatrix}
0&0&0\\
0&0&-1\\
0&1&0
\end{pmatrix},
y
\begin{pmatrix}
0&0&0&\cdots&0\\
0&0&-X_{2}&\cdots&-X_{n}\\
0&X_{2}&0&\cdots&0\\
\vdots&\vdots&\vdots&\cdots&\vdots\\
0&X_{n}&0&\cdots&0
\end{pmatrix}
\right)
\end{equation*}
if and only if $yX_{2}=\cdots=yX_{n}=0$. Suppose that
$[{v_{1}},{v_{2}}]\in{\frak c}({\frak k})$. Then
\begin{equation*}
\Vert[{v_{1}},{v_{2}}]\Vert^{2}
=(xy^{\prime})^{2}+\sum^{n}_{i=2}(yX_{i})^{2}
=(xy^{\prime})^{2}\leq{1}.
\end{equation*}
Here $\Vert[{v_{1}},{v_{2}}]\Vert^{2}=1$ if and only if
\begin{equation*}
x^{2}=1 \text{ and } (y^{\prime})^{2}=1, \text{ i.e., } x=1 \text{
and } y^{\prime}=\pm{1},
\end{equation*}
if and only if
\begin{equation*}
v_{1}=\left(
\begin{pmatrix}
0&1&0\\
-1&0&0\\
0&0&0
\end{pmatrix},
0\right)=H_{2},\quad v_{2}=\left(
\begin{pmatrix}
0&0&\pm{1}\\
0&0&0\\
\mp{1}&0&0
\end{pmatrix},
0\right)=\pm{JH_{2}}.
\end{equation*}
This implies the statement (c).

(d) Let $(U,K)=(SO(m+2),SO(2)\times SO(m)){\ }(n=2m-2,m\geq{3})$.
Suppose that $\mu([V])\in{\frak c}({\frak k})$. We can choose an
orthonormal basis $\{v_{1},v_{2}\}$ of $V$ compatible with the
orientation of $[V]$ such that $v_{1}\in{\frak a}\cap{S^{n+1}(1)}$
and $(-\sqrt{-1})\gamma(v_{1})\geq{0}$ for each
$\gamma\in\Sigma^{+}(U,K)$. We express $v_{2}$ as
\begin{equation*}
v_{2}=v_{2,0}+\sum_{\gamma\in{\Sigma^{+}(U,K)}}v_{2,\gamma}
\end{equation*}
where $v_{2,0}\in{\frak a}$ and $v_{2,\gamma}\in{\frak p}_{\gamma}$.
In this case
\begin{equation*}
\begin{split}
\Sigma^{+}(U,K) =\{{\
}&\gamma_{1}=\sqrt{-1}(\varepsilon_{1}-\varepsilon_{2}), {\
}\gamma_{2}=\sqrt{-1}\varepsilon_{2},
{\ }\gamma_{1}+\gamma_{2}=\sqrt{-1}\varepsilon_{1}, \\
&\tilde{\gamma}=\gamma_{1}+2\gamma_{2}
=\sqrt{-1}(\varepsilon_{1}+\varepsilon_{2}){\ }\}
\end{split}
\end{equation*}
and set

$\{H_{\gamma_{1}}=e_{1}-e_{2}, H_{\gamma_{2}}=e_{2},
H_{\gamma_{1}+\gamma_{2}}=e_{1},
H_{\tilde{\gamma}}=H_{\gamma_{1}+2\gamma_{2}}
=e_{1}+e_{2}\}\subset{\frak a}$.

Then $\{H_{1}:=e_{2},H_{2}:=e_{1}\}$ is an orthonormal basis of
${\frak a}$, which defines an orientation of ${\frak a}$, and we
have the root decomposition of ${\frak a}^{\perp}$ with respect to
$\Sigma^{+}(U,K)$ as
\begin{equation*}
{\frak a}^{\perp}=\sum_{\gamma\in\Sigma^{+}(U,K)}{\frak p}_{\gamma}
={\frak p}_{\gamma_{1}}+{\frak p}_{\gamma_{2}} +{\frak
p}_{\gamma_{1}+\gamma_{2}}+{\frak p}_{\gamma_{1}+2\gamma_{2}}{\ },
\end{equation*}
where $\dim{\frak p}_{\gamma_{1}}=1$, $\dim{\frak
p}_{\gamma_{2}}=m-2$, $\dim{\frak p}_{\gamma_{1}+\gamma_{2}}=m-2$,
$\dim{\frak p}_{\gamma_{1}+2\gamma_{2}}=1$. Using the orthonormal
basis \eqref{StdBasis1},\eqref{StdBasis2}, we express each
$v_{2,\gamma}$ as
\begin{equation*}
v_{2,\gamma}=\sum^{m(\gamma)}_{i=1}c_{\gamma,i}Y_{\gamma,i}{\ },
\end{equation*}
where $c_{\gamma,i}\in{\bold R}$. Then we have
\begin{equation*}
v_{2}=v_{2,0}+c_{\gamma_{1},1}Y_{\gamma_{1},1}
+c_{\tilde{\gamma},1}Y_{\tilde{\gamma},1}
+\sum^{m-2}_{i=1}c_{\gamma_{2},i}Y_{\gamma_{2},i}
+\sum^{m-2}_{i=1}c_{\gamma_{1}+\gamma_{2},i}
Y_{\gamma_{1}+\gamma_{2},i}.
\end{equation*}
and
\begin{equation*}
\Vert{v_{2}}\Vert^{2}
=\Vert{v_{2,0}}\Vert^{2}+\vert{c_{\gamma_{1},1}}\vert^{2}
+\vert{c_{\tilde{\gamma},1}}\vert^{2}
+\sum^{m-2}_{i=1}\vert{c_{\gamma_{2},i}}\vert^{2}
+\sum^{m-2}_{i=1}\vert{c_{\gamma_{1}+\gamma_{2},i}}\vert^{2} =1.
\end{equation*}
Thus
\begin{equation*}
\begin{split}
&[v_{1},v_{2}]\in{\frak c}({\frak k}) ={\bold
R}(X_{\gamma_{1}}+X_{\tilde{\gamma}})
\\
=& c_{\gamma_{1},1}[v_{1},Y_{\gamma_{1},1}]
+c_{\tilde{\gamma},1}[v_{1},Y_{\tilde{\gamma},1}]
+\sum^{m-2}_{i=1}c_{\gamma_{2},i}[v_{1},Y_{\gamma_{2},i}]
+\sum^{m-2}_{i=1}c_{\gamma_{1}+\gamma_{2},i}
[v_{1},Y_{\gamma_{1}+\gamma_{2},i}]\\
=& c_{\gamma_{1},1}(-\sqrt{-1}\gamma_{1}(v_{1})X_{\gamma_{1}})
+c_{\tilde{\gamma},1}(-\sqrt{-1}\tilde{\gamma}(v_{1})X_{\tilde{\gamma}})\\
=& -c_{\gamma_{1},1}(\sqrt{-1}\gamma_{1}(v_{1}))X_{\gamma_{1}}
-c_{\tilde{\gamma},1}(\sqrt{-1}\tilde{\gamma}(v_{1}))X_{\tilde{\gamma}}\\
=& -c_{\gamma_{1},1}(v_{1,1}-v_{1,2})X_{\gamma_{1}}
-c_{\tilde{\gamma},1}(v_{1,1}+v_{1,2})X_{\tilde{\gamma}},
\end{split}
\end{equation*}
where $v_{1}=v_{1,1}e_{1}+v_{1,2}e_{2} \in{\frak a}$. Note that
$v_{1,1}^{2}+v_{1,2}^{2}=1$, $v_{1,1}\geq{0}$, $v_{1,2}\geq{0}$,
$v_{1,1}-v_{1,2}\geq{0}$, and $v_{1,1}+v_{1,2}\geq{0}$. Hence we
have
\begin{equation*}
c_{\gamma_{1}}(v_{1,1}-v_{1,2}) =c_{\tilde{\gamma}}(v_{1,1}+v_{1,2})
\end{equation*}
and
\begin{equation*}
\vert{c_{\gamma_{1},1}}\vert(v_{1,1}-v_{1,2})
=\vert{c_{\tilde{\gamma},1}}\vert(v_{1,1}+v_{1,2}).
\end{equation*}
It implies that
$\vert{c_{\gamma_{1},1}}\vert\geq\vert{c_{\tilde{\gamma},1}}\vert$.
Therefore we obtain
\begin{equation}
\begin{split}
\Vert{[v_{1},v_{2}]}\Vert^{2} =& \vert{c_{\gamma_{1}}}\vert^{2}\vert
v_{1,1}-v_{1,2}\vert^{2}
+\vert{c_{\tilde{\gamma}}}\vert^{2}\vert v_{1,1}+v_{1,2}\vert^{2}\\
=& \vert{c_{\gamma_{1},1}}\vert^{2}
(v_{1,1}^{2}-2v_{1,1}v_{1,2}+v_{1,2}^{2})\\
&+\vert{c_{\tilde{\gamma},1}}\vert^{2}
(v_{1,1}^{2}+2v_{1,1}v_{1,2}+v_{1,2}^{2})\\
=& \vert{c_{\gamma_{1},1}}\vert^{2} (1-2v_{1,1}v_{1,2})
+\vert{c_{\tilde{\gamma},1}}\vert^{2}
(1+2v_{1,1}v_{1,2})\\
=& \vert{c_{\gamma_{1},1}}\vert^{2}
+\vert{c_{\tilde{\gamma},1}}\vert^{2} +2v_{1,1}v_{1,2}
(\vert{c_{\tilde{\gamma},1}}\vert^{2}
-\vert{c_{\gamma_{1},1}}\vert^{2})\\
\leq& 1+2v_{1,1}v_{1,2} (\vert{c_{\tilde{\gamma},1}}\vert^{2}
-\vert{c_{\gamma_{1},1}}\vert^{2})\leq{1}.
\end{split}
\end{equation}
Here $\Vert{[v_{1},v_{2}]}\Vert^{2}=1$ if and only if
$v_{1}=e_{1}=H_{2}$ and $v_{2}=\pm\frac{1}{\sqrt{2}}
(Y_{\gamma_{1},1}+Y_{\tilde{\gamma},1})={\pm}JH_{2}$. This implies
the statement (d).
\end{proof}

\subsection{Two one-parameter families of Lagrangian orbits
in complex hyperquadrics}\label{subsec:4.2}

In the cases of $(U,K)=(S^{1}\times SO(3),SO(2))$ and
$(U,K)=(SO(3)\times SO(3),SO(2)\times SO(2))$, the induced action of
$K=SO(2)$ on $Q_{1}({\bold C})$ is just the standard $S^{1}$-action
on $S^{2}$ via the identification $Q_{1}({\bold C})\cong S^{2}$ and
the induced action of $K=SO(2)\times SO(2)$ on $Q_{2}({\bold C})$ is
the product standard $S^{1}\times S^{1}$-action on $S^{2}\times
S^{2}$ via the identification $Q_{2}({\bold C})\cong S^{2}\times
S^{2}$, respectively. In the first case every $K$-orbit except for
two fixed points ($0$-dimensional isotropic orbits) is a Lagrangian
orbit. In the second case every $K$-orbit except for four fixed
points ($0$-dimensional isotropic orbits) and for four one-parameter
families of $1$-dimensional isotropic orbits is a Lagrangian orbit.

Suppose that $(U,K)$ is
$(SO(3)\times SO(n+1),SO(2)\times SO(n)){\ }(n\geq{3})$ or
$(SO(m+2),SO(2)\times SO(m)){\ }(n=2m-2,m\geq{3})$.

Let $\{H_{1},H_{2}\}$ be an orthonormal basis of ${\frak a}$, which
defines an orientation of ${{\frak a}}$, described in \eqref{H1H2}
or the proof of Lemma \ref{NormMomentMap} (c),(d). For each
$\lambda=e^{\sqrt{-1}\theta}\in S^{1}$, we define a $2$-dimensional
vector subspace $W_{\lambda}$ of ${\frak p}$ by
\begin{equation}
W_{\lambda}:={\bold R}(\cos\theta H_{1}+\sin\theta JH_{2}) +{\bold
R}H_{2}.
\end{equation}
We denote by $[W_{\lambda}]\in\widetilde{\mathrm {Gr}}_{2}({\frak
p})$ a $2$-dimensional vector subspace $W_{\lambda}$ oriented by the
basis $\{\cos\theta H_{1}+\sin\theta JH_{2},H_{2}\}$. Then by Lemma
\ref{[JH_{2},H_{2}]}
we obtain
\begin{equation}
\begin{split}
\mu([W_{\lambda}])
&=-[\cos\theta H_{1}+\sin\theta JH_{2},H_{2}]\\
&=-\sin\theta [JH_{2},H_{2}]\in{\frak c}({\frak k}).
\end{split}
\end{equation}
Hence for each $\lambda\in{S^{1}}$ the $K$-orbit
$K\cdot[W_{\lambda}]$ in $\widetilde{\mathrm {Gr}}_{2}({\frak
p})=Q_{n}({\bold C})$ is an isotropic orbit. More precisely,
$K\cdot[W_{\lambda}]$ is a Lagrangian orbit if
$\lambda\not=\pm\sqrt{-1}$, and an isotropic orbit of dimension less
than $n$ if $\lambda=\pm\sqrt{-1}$. If $(U,K)=(SO(3)\times
SO(n+1),SO(2)\times SO(n)){\ }(n\geq{2})$, then
$\dim(K\cdot[W_{\pm\sqrt{-1}}])=0$. If $(U,K)=(SO(m+2),SO(2)\times
SO(m)){\ }(n=2m-2,m\geq{3})$, then
\begin{enumerate}
\renewcommand{\labelenumi}{(\arabic{enumi})}
\item
$\dim(K\cdot[W_{\pm\sqrt{-1}}])=m-1$ and $K\cdot[W_{\pm\sqrt{-1}}]$
is an isotropic orbit diffeomorphic to $SO(m)/S(O(1)\times O(m-1))
\cong{\bold R}P^{m-1}$,
\item
for each $\lambda\in{S^{1}\setminus\{\pm{\sqrt{-1}}\}}$,
$\dim(K\cdot[W_{\lambda}])=n=2m-2$ and
$K\cdot[W_{\lambda}]=K\cdot[W_{-\bar{\lambda}}]$ is a Lagrangian
orbit diffeomorphic to
$(SO(2)\times{SO(m)})/({\bold Z}_{2}\times{\bold Z}_{4}\times{SO(m-2)})
\cong N^{n}/{\bold Z}_{4}$.
\end{enumerate}

\subsection{Classification theorem of Lagrangian orbits
in complex hyperquadrics}\label{subsec:4.3}

In the case (a) or (b), we already know all Lagrangian orbits. In
the case of (c) or (d), as we showed in
Subsection \ref{subsec:4.2}, there exists an $S^{1}$-family
$[W_{\lambda}]\in\widetilde{\mathrm {Gr}}_{2}({\frak p})$ with
$\lambda\in{S^{1}=\{\lambda\in{\bold C}{\ }\vert{\ }
\vert{\lambda}\vert=1\}}$ satisfying the following conditions :
\begin{enumerate}
\renewcommand{\labelenumi}{(\arabic{enumi})}
\item
$[W_{1}]=[{\frak a}]$, $[W_{-\lambda}]=-[W_{\lambda}]$.
\item
$\mu([W_{\lambda}])\in{\frak c}({\frak k})$ and
$\mu([W_{-\lambda}])=-\mu([W_{\lambda}])$.
\item
For each $\lambda\in{S^{1}}$ with $\lambda\not=\pm\sqrt{-1}$, the
$K$-orbit $K\cdot{[W_{\lambda}]}$ is a Lagrangian orbit in
$Q_{n}({\bold C})$, and the $K$-orbits $K\cdot{[W_{\pm\sqrt{-1}}]}$
are isotropic orbits in $Q_{n}({\bold C})$ with
$\dim{K\cdot{[W_{\pm\sqrt{-1}}]}}=0$ in case (c) and
$\dim{K\cdot{[W_{\pm\sqrt{-1}}]}}=m-1$ in case (d).
\item
The square norm $\Vert{\mu}\Vert^{2}$ of the moment map satisfies
\begin{equation}
\begin{split}
&\mathrm{Max}\{\Vert{\mu([W_{\lambda}])}\Vert^{2}{\ }\vert{\ }
\lambda\in{S^{1}}\}\\
=&\Vert{\mu([W_{\sqrt{-1}}])}\Vert^{2}=\Vert{\mu([W_{-\sqrt{-1}}])}\Vert^{2}\\
=&\mathrm{Max}\{\Vert{\mu([V])}\Vert^{2}{\ }\vert{\ }
[V]\in\widetilde{\mathrm Gr}_{2}({\frak p}) \text{ with
}\mu([V])\in{\frak c}({\frak k})\}.
\end{split}
\end{equation}
\end{enumerate}
Thus we have
\begin{equation*}
\Vert{\mu([W_{\pm{1}}])}\Vert^{2} =\Vert{\mu(\pm[{\frak
a}])}\Vert^{2} =0 \leq \Vert{\mu([V_{0}])}\Vert^{2} \leq
\Vert{\mu([W_{\pm\sqrt{-1}}])}\Vert^{2}.
\end{equation*}
Since $\mathrm{Im}\mu\cap{\frak c}({\frak k})$ is a $1$-dimensional
connected compact subset by Proposition \ref{KirwanResult} (2),
there is $\lambda_{0}\in{S^{1}}$ such that
$\Vert{\mu([V_{0}])}\Vert^{2}
=\Vert{\mu([W_{\lambda_{0}}])}\Vert^{2}$. Therefore we obtain
$\xi=\mu([V_{0}])=\mu([W_{\lambda_{0}}])$ or
$\xi=\mu([V_{0}])=\mu([W_{-\lambda_{0}}]$. Since $\mu^{-1}(\xi)$ is
connected by Proposition \ref{KirwanResult} (1), we conclude that
$K\cdot{[V_{0}]}=\mu^{-1}(\xi)=K\cdot{[W_{\lambda_{0}}]}$ or
$K\cdot{[V_{0}]}=\mu^{-1}(\xi)=K\cdot{[W_{-\lambda_{0}}]}$.

Therefore our classification of compact homogeneous Lagrangian
submanifolds in complex hyperquadrics are described as follows.

\begin{theorem}
Let $L$ be a compact homogeneous Lagrangian submanifold in
$Q_{n}({\bold C})$. Then there exists uniquely a homogeneous
isoparametric hypersurface $N^{n}$ in $S^{n+1}(1)$ obtained as a
principal orbit of the isotropy action of a compact Riemannian
symmetric pair $(U,K)$ of rank $2$ such that the following
statements hold :
\begin{enumerate}
\renewcommand{\labelenumi}{(\arabic{enumi})}
\item
If $(U,K)$ is not one of
\begin{enumerate}
\renewcommand{\labelenumi}{(\roman{enumi})}
\item $(S^{1}\times SO(3),SO(2))$,
\item $(SO(3)\times SO(3),SO(2)\times SO(2))$,
\item $(SO(3)\times SO(n+1),SO(2)\times SO(n)){\ }(n\geq{3})$,
\item $(SO(m+2),SO(2)\times SO(m)){\ }(n=2m-2,m\geq{3})$,
\end{enumerate}
then ${\frak c}({\frak k})\cap\mathrm{Im}\mu=\{0\}$ and
\begin{equation*}
L={\mathcal G}(N^{n})\subset{Q_{n}({\bold C})},
\end{equation*}
which is a minimal Lagrangian submanifold in $Q_{n}({\bold C})$.
\item
If $(U,K)$ is $(S^{1}\times SO(3),SO(2))$, then $L$ is a small or
great circle in $Q_{1}({\bold C})\cong{S^{2}}$.
\item
If $(U,K)$ is $(SO(3)\times SO(3),SO(2)\times SO(2))$, then $L$ is a
product of small or great circles of $S^{2}$ in $Q_{2}({\bold
C})\cong{S^{2}\times S^{2}}$.
\item
If $(U,K)$ is $(SO(3)\times SO(n+1),SO(2)\times SO(n))
{\ }(n\geq{2})$ , then
$$
\mathrm{Im}\mu\cap{\frak c}({\frak k})=
\{\mu([W_{\lambda}]){\ }\vert{\ }\lambda\in{S^{1}}\}
$$
and
\begin{equation*}
L=K\cdot{[W_{\lambda}]}\subset{Q_{n}({\bold C})}
\quad\text{ for some }\lambda\in{S^{1}\setminus\{\pm{\sqrt{-1}}\}},
\end{equation*}
where $K\cdot{[W_{\lambda}]}{\ }(\lambda\in{S^{1}})$ is the
$S^{1}$-family of Lagrangian or isotropic $K$-orbits satisfying
\begin{enumerate}
\renewcommand{\labelenumi}{(\alph{enumi})}
\item
$K\cdot{[W_{1}]}=K\cdot{[W_{-1}]}={\mathcal G}(N^{n})$ is a totally
geodesic Lagrangian submanifold in $Q_{n}({\bold C})$.
\item
For each $\lambda\in{S^{1}\setminus\{\pm{\sqrt{-1}}\}}$,
$K\cdot{[W_{\lambda}]}$ is a Lagrangian orbit in $Q_{n}({\bold C})$
which is diffeomorphic to $(S^{1}\times{S^{n-1}})/{\bold Z}_{2}
\cong{Q_{2,n}({\bold R})}$.
\item
$K\cdot{[W_{\pm\sqrt{-1}}]}$ are isotropic orbits in $Q_{n}({\bold
C})$ with $\dim{K\cdot{[W_{\pm\sqrt{-1}}]}}=0$.
\end{enumerate}
\item
If $(U,K)$ is $(SO(m+2),SO(2)\times SO(m)){\ }(n=2m-2)$, then
$$
\mathrm{Im}\mu\cap{\frak c}({\frak k})= \{\mu([W_{\lambda}]){\
}\vert{\ }\lambda\in{S^{1}}\}
$$
and
\begin{equation*}
L=K\cdot{[W_{\lambda}]}\subset{Q_{n}({\bold C})}
\quad\text{ for some }\lambda\in{S^{1}\setminus\{\pm{\sqrt{-1}}\}},
\end{equation*}
where $K\cdot{[W_{\lambda}]}{\ }(\lambda\in{S^{1}})$ is the
$S^{1}$-family of Lagrangian or isotropic orbits satisfying
\begin{enumerate}
\renewcommand{\labelenumi}{(\alph{enumi})}
\item
$K\cdot{[W_{1}]}=K\cdot{[W_{-1}]}={\mathcal G}(N^{n})$ is a minimal
Lagrangian submanifold in $Q_{n}({\bold C})$.
\item
For each $\lambda\in{S^{1}\setminus\{\pm{\sqrt{-1}}\}}$,
$K\cdot{[W_{\lambda}]}$ is a Lagrangian orbit in $Q_{n}({\bold C})$,
which is diffeomorphic to
$(SO(2)\times SO(m))/({\bold Z}_{2}\times{\bold Z}_{4}\times{SO(m-2)})$.
\item
$K\cdot{[W_{\pm\sqrt{-1}}]}$ are isotropic orbits in $Q_{n}({\bold
C})$ with $\dim{K\cdot{[W_{\pm\sqrt{-1}}]}}=m-1$ which is
diffeomorphic to $SO(m)/S(O(1)\times O(m-1))\cong{\bold R}P^{m-1}$.
\end{enumerate}
\end{enumerate}
\end{theorem}

\begin{rem0}
In each case when
$(U,K)=(SO(3)\times SO(n+1),SO(2)\times SO(n)){\ }(n\geq{3})$
or $(U,K)=(SO(m+2),SO(2)\times SO(m))$ with $n=2m-2$,
there is a nontrivial one-parameter family of Lagrangian orbits in
$Q_{n}({\bold C})$.
The family contains homogeneous Lagrangian
submanifolds which can NEVER be obtained as the Gauss images of
homogeneous isoparametric hypersurfaces.
\end{rem0}

\begin{corollary}
Any compact homogeneous Lagrangian submanifold in a complex
hyperquadric is obtained as the Gauss image of a compact homogeneous
isoparametric hypersurface in a sphere, or as its Lagrangian
deformation.
\end{corollary}

\section{Hamiltonian Stability of Gauss images of isoparametric
hypersurfaces
in spheres}\label{sec:5}

Let $N^{n}$ be an oriented compact isoparametric hypersurface
embedded in $S^{n+1}(1)$. Now we already know that its Gauss map
${\mathcal G}:N^{n}\rightarrow Q_{n}({\bold C})$ is a minimal
Lagrangian immersion. In \cite{Palmer97} Palmer showed that the
Gauss map ${\mathcal G}:N^{n}\rightarrow Q_{n}({\bold C})$ is
Hamiltonian stable if and only if $N^{n}=S^{n}\subset S^{n+1}(1)$
($g=1$).
\begin{prob1}
Investigate the Hamiltonian stability of its Gauss image ${\mathcal
G}(N^{n})=N^{n}/{\bold Z}_{g}$ embedded in ${Q_{n}({\bold C})}$ as a
compact minimal Lagrangian submanifold.
\end{prob1}

Let $g^{\mathrm{std}}_{Q_{n}({\bold C})}$ denote the
$SO(n+2)$-invariant Riemannian metric induced from the standard
Euclidean metric of ${\bold R}^{n+2}$, whose Einstein constant is
equal to $n$. Let $g^{\mathrm{KC}}_{Q_{n}({\bold C})}$ denote the
$SO(n+2)$-invariant Riemannian metric induced from the
Killing-Cartan form of $SO(n+2)$, whose Einstein constant is equal
to $\frac{1}{2}$ (cf.~\cite{KNI-II}). We also can use
$\mathrm{Table}{\ }1$ in Section \ref{sec:1} to determine the
Hamiltonian stability in the cases of $g=1$ and $g=2$, because
${\mathcal G}(N^{n})$ is a totally geodesic Lagrangian submanifold
in $Q_{n}({\bold C})$ in these cases.

\noindent $g=1$ : ${\mathcal G}(N^{n})=Q_{1,n+1}({\bold
R})\subset{Q_{n}({\bold C})}$ is Hamiltonian stable.

\noindent $g=2$ : $N^{n}=S^{m_{1}}\times S^{m_{2}} {\ }(1\leq
m_{1}\leq m_{2})$ are the so called Clifford hypersurfaces.

If $m_{2}-m_{1}\geq{3}$, then ${\mathcal
G}(N^{n})=Q_{m_{1}+1,m_{2}+1}({\bold R}) \subset{Q_{n}({\bold C})}$
is Not Hamiltonian stable. Otherwise ${\mathcal
G}(N^{n})=Q_{m_{1}+1,m_{2}+1}({\bold R}) \subset{Q_{n}({\bold C})}$
is Hamiltonian stable.

In the case of $g=3$, all isoparametric hypersurfaces are
homogeneous by E.~Cartan's result. In this section we prove the
following result.

\begin{theorem}
Suppose that $g=3$, that is, $N^{n}$ is one of the following
isoparametric hypersurfaces :
$(1){\ }SO(3)/{\bold Z}_{2}+{\bold Z}_{2}$,
$(2){\ }SU(3)/T^{2}$,
$(3){\ }Sp(3)/Sp(1)^{3}$,
$(4){\ }F_{4}/Spin(8)$.
Then $L={\mathcal G}(N^{n})\subset{Q_{n}({\bold C})}$ is
strictly Hamiltonian stable.
\end{theorem}

\begin{rem0}
In case $g=3$ the induced metrics from $Q_{n}({\bold C})$ have nice
intrinsic properties. In (1), $L$ has constant sectional curvature
$1/96$ if the Einstein constant of $Q_{n}({\bold C})$ is equal to
$\frac{1}{2}$. In (1)--(4), $L$ has non-negative sectional
curvatures. In fact, the induced invariant metrics are normal
homogenous metrics.
\end{rem0}

Refer also \cite{MOU84}, \cite{Kotani85}, \cite{Muto88} for
investigation on the first eigenvalue of the Laplacian of
homogeneous isoparametric hypersurfaces in spheres.

We recall some results from the spherical function theory on compact
homogeneous spaces to order to determine the first eigenvalues of
${\mathcal G}(N^{n})$ relative to the induced metric from
$Q_{n}({\bold C})$,

Suppose that $\langle{\ },{\ }\rangle_{\frak k}$ is an
$\mathrm{Ad}K$-invariant inner product of ${\frak k}$. For a compact
Lie subgroup $S$ of $K$ with Lie algebra ${\frak s}$, we take the
orthogonal direct sum decomposition ${\frak k}={\frak s}+{\frak m}$
and the vector space ${\frak m}$ is identified with the tangent
vector space $T_{eS}(K/S)$ of compact homogeneous space $K/S$ at the
origin $eS$. The {\it Casimir operator} ${\mathcal C}$ of $(K,S)$
with respect to $\langle{\ },{\ }\rangle_{\frak k}$ by
${\mathcal C}:=\sum^{n}_{i=1}X_{i}$,
where $\{X_{i}{\ }\vert{\ }i=1,\cdots,n\}$ is an orthonormal basis
of ${\frak m}$ with respect to $\langle{\ },{\ }\rangle_{\frak k}$.
Let ${\mathcal D}(K)$ be the complete set of all inequivalent
irreducible unitary representations of a compact Lie group $K$. For
a maximal abelian subalgebra ${\frak h}$ of ${\frak k}$, let
$\Sigma(K)$ be the set of all roots of ${\frak k}$ with respect to
${\frak h}$ and let $\Sigma^{+}(K)$ be its subset of all positive
$\alpha\in\Sigma(K)$ relative to a linear order on ${\frak h}$. Set
\begin{equation*}
\begin{split}
\Gamma(K)&:=\{\xi\in{\frak h}{\ }\vert{\ }\exp(\xi)=e\},\\
Z(K)&:=\{\Lambda\in{\frak h}^{\ast}{\ }\vert{\ }\Lambda(\xi)\in{\bold
Z}\},\\
D(K)&:=\{\Lambda\in{\frak h}^{\ast}{\ }\vert{\
}(\Lambda,\alpha)\geq{0} \text{ for each
}\alpha\in{\Sigma^{+}(K)}\}.
\end{split}
\end{equation*}
Then we know that there is a bijective correspondence between $D(K)$
and ${\mathcal D}(K)$ : Each $\Lambda\in{D(K)}$ uniquely corresponds
to an irreducible unitary representation
$(V_{\Lambda},\rho_{\Lambda})$ of $K$ with the highest weight
$\Lambda$, up to the equivalence. Here we denote by $\langle{\ },{\
}\rangle_{V_{\Lambda}}$ a $K$-invariant Hermitian inner product
equipped on $V_{\Lambda}$. Let $(V_{\Lambda})_{S}$ denote the vector
subspace of $V_{\Lambda}$ consisting of all vectors fixed by
$\rho_{\Lambda}(S)$. Define
\begin{equation*}
D(K,S):=\{\Lambda\in{D(K)}{\ }\vert{\ }
(V_{\Lambda})_{S}\not=\{0\}\}.
\end{equation*}
By the Peter-Weyl's theorem we know that
\begin{equation*}
C^{\infty}(K/S) =\bigoplus_{\Lambda\in{D(K,S)}}
(V_{\Lambda})^{\ast}_{S}\otimes{V_{\Lambda}}.
\end{equation*}
Here each $w\in(V_{\Lambda})_{S}$ and each $v\in{V_{\Lambda}}$
correspond to $f\in{C^{\infty}(K/S)}$ defined by
\begin{equation*}
f(aS):=\langle{\rho_{\Lambda}(a)w},{v}\rangle_{\Lambda}
\quad(aK\in{K/S}).
\end{equation*}
Then the Laplace-Beltrami operator $\Delta_{K/S}$ with respect to
the metric $g_{K/S}$ on $K/S$ induced by $\langle{\ },{\
}\rangle_{\frak k}$ is expressed in terms of ${\mathcal C}$ as
\begin{equation*}
(\Delta_{K/S}{f})(aS) =\langle {\rho_{\Lambda}(a)
\left((d\rho_{\Lambda}({\mathcal C})w)\right)},{v}
\rangle_{\Lambda}.
\end{equation*}
By Schur's lemma there is a real constant $c(\Lambda,\langle{\ },{\
}\rangle_{\frak k})\leq{0}$ such that
\begin{equation*}
(d\rho_{\Lambda}({\mathcal C}))v =\sum^{n}_{i=1}
(d\rho_{\Lambda}(X_{\gamma,i}))^{2}v =c(\Lambda,\langle{\ },{\
}\rangle_{\frak k})v \qquad\text{ for each } v\in{V_{\lambda}}.
\end{equation*}
The eigenvalue $c(\Lambda,\langle{\ },{\ }\rangle_{\frak k})$ is
described by the Freudenthal's formula
\begin{equation*}
c(\Lambda,\langle{\ },{\ }\rangle_{\frak k})
=-\langle{\Lambda},{\Lambda+2\delta}\rangle_{\frak k},
\end{equation*}
where $2\delta=\sum_{\alpha\in{\Sigma^{+}(K)}}\alpha$.

We shall consider our compact homogeneous spaces $K/K_{0}$ and
$K/K_{[{\frak a}]}$. The tangent vector spaces $T_{eK_{0}}(K/K_{0})$
and $T_{eK_{[\frak a]}}(K/K_{[{\frak a}]})$ can be identified with
the vector subspace ${\frak m}$ in \eqref{M&Aperp} of ${\frak k}$.
Let $\langle{\ },{\ }\rangle$ denote the $K_{\frak a}$-invariant
inner product of ${\frak m}$ corresponding to the $K$-invariant
Riemannian metric ${\mathcal G}^{\ast}g^{std}_{Q_{n}({\bold C})}$ on
$K/K_{0}$ induced from $g^{std}_{Q_{n}({\bold C})}$ through the
Gauss map ${\mathcal G}$. Using the standard basis
\eqref{StdBasis1},\eqref{StdBasis2}, we compute
\begin{equation*}
\begin{split}
[(d{\mathcal G})_{eK_{0}}(X_{\gamma,i})](H) =\pi_{{\frak
a}^{\perp}}([X_{\gamma,i},H])
=[X_{\gamma,i},H]=-\sqrt{-1}\gamma(H)Y_{\gamma,i}
\end{split}
\end{equation*}
and
\begin{equation*}
\begin{split}
\langle (d{\mathcal G})_{eK_{0}}(X_{\gamma,i}), (d{\mathcal
G})_{eK_{0}}(X_{\gamma^{\prime},j}) \rangle_{\frak u}
=\Vert{\gamma}\Vert_{\frak u}^{2} \langle
X_{\gamma,i},X_{\gamma^{\prime},j} \rangle_{\frak u}{\ }.
\end{split}
\end{equation*}
Thus we see that
\begin{equation*}
\{{\ }\frac{1}{\Vert{\gamma}\Vert_{\frak u}}X_{\gamma,i} {\ }\vert{\
}\gamma\in{\Sigma^{+}(U,K)},i=1,2,\cdots,m(\gamma){\ }\}
\end{equation*}
is an orthonormal basis of ${\frak m}$ with respect to
$\langle{\ },{\ }\rangle$.

Suppose that $g=3$. Since $(U,K)$ is type $A_{2}$, we have
\begin{equation}
\langle{\ },{\ }\rangle=\Vert{\gamma_{1}}\Vert_{\frak u}^{2}
\langle{\ },{\ }\rangle_{\frak u}
\end{equation}
on ${\frak m}$. On the other hand, since $K$ is simple, we can
define an $\mathrm{Ad}K$-invariant inner product $\langle{\ },{\
}\rangle_{\frak k}=-B_{\frak k}({\ },{\ })$ of ${\frak k}$ by using
the Killing-Cartan form $B_{\frak k}$ of ${\frak k}$ and moreover
there is $b>0$ such that
\begin{equation}
\langle{\ },{\ }\rangle_{\frak k} =b\langle{\ },{\ }\rangle_{\frak
u}.
\end{equation}
on ${\frak k}$.
Set $C:=\Vert{\gamma_{1}}\Vert_{\frak u}^{2}\cdot{b^{-1}}$.
Then we have
\begin{equation}
\langle{\ },{\ }\rangle= C\langle{\ },{\ }\rangle_{\frak k}.
\end{equation}
Thus we obtain

\begin{lemma}
$L={\mathcal G}(N^{n})\subset{Q_{n}({\bold C})}$ is Hamiltonian
stable if and only if
\begin{equation*}
\mathrm{Min} \{-c(\Lambda,\langle{\ },{\ }\rangle_{\frak k}) {\
}\vert{\ }0\not=\Lambda\in{D(K,K_{[\frak a]})}\}
\end{equation*}
is equal to $Cn$.
\end{lemma}

\begin{lemma}
The constants $b^{-1}$, $\Vert{\gamma_{1}}\Vert_{\frak u}^{2}$, $Cn$
in each case are given as in the following table :

\begin{center}
\begin{tabular}
{|c|c|c|c|c|c|c|} \hline $(U,K)$ &$n$ &$\dim{\frak p}$ &
$\dim{\frak k}$
&$b^{-1}$ &$\Vert{\gamma_{1}}\Vert_{\frak u}^{2}$ &$Cn$
\\
\hline $(SU(3),SO(3))$ &$3$ &$5$ &$3$
&$6$ &$1/3$ &$6$
\\
\hline $(SU(3)\times SU(3),SU(3))$ &$6$ &$8$ &$8$
&$2$ &$1/6$ &$2$
\\
\hline $(SU(6),Sp(3))$ &$12$ &$14$ &$21$
&$3/2$ &$1/12$ &$3/2$
\\
\hline $(E_{6},F_{4})$ &$24$ &$26$ &$52$
&$4/3$ &$1/24$ &$4/3$
\\
\hline
\end{tabular}
\end{center}

\end{lemma}

\begin{proof}
Choose an orthonormal basis $\{e_{j}\}$ of ${\frak m}$ and an
orthonormal basis $\{H_{\nu}{\ }\vert{\ }\nu=1,2\}$ of ${\frak a}$
with respect to $\langle{\ },{\ }\rangle_{\frak u}$. Then we compute
\begin{equation*}
\sum_{j}\langle{e_{j}},{e_{j}}\rangle
=Cb\sum_{j}\langle{e_{j}},{e_{j}}\rangle_{\frak u}=Cb\dim{\frak m},
\end{equation*}
\begin{equation*}
\begin{split}
\sum_{j}\langle{e_{j}},{e_{j}}\rangle
=&\sum_{\nu}\sum_{j}\langle{[e_{j},H_{\nu}]},
{[e_{j},H_{\nu}]}\rangle_{\frak u}
=- \sum_{\nu}\langle\sum_{j}(\mathrm{ad}e_{j})^{2}H_{\nu},
H_{\nu}\rangle_{\frak u}\\
=&-b
\sum_{\nu}\langle\sum_{j}(\mathrm{ad}e_{j}/\sqrt{b})^{2}H_{\nu},
H_{\nu}\rangle_{\frak u} =-2b\cdot c({\frak k},\mathrm{ad}_{\frak
p},\langle{\ },{\ }\rangle_{\frak k}),
\end{split}
\end{equation*}
and
\begin{equation*}
\begin{split}
\sum_{j}\langle{e_{j}},{e_{j}}\rangle =&
-\sum_{\nu}\sum_{j}\langle{(\mathrm{ad}H_{\nu})^{2}e_{j}},
{e_{j}}\rangle_{\frak u}
=-\sum_{\nu}\mathrm{tr}_{\frak k}(\mathrm{ad}H_{\nu})^{2}\\
=&-\frac{1}{2}\sum_{\nu}\mathrm{tr}_{\frak
u}(\mathrm{ad}H_{\nu})^{2}
=\frac{1}{2}\sum_{\nu}\langle{H_{\nu}},{H_{\nu}}\rangle_{\frak u}
=1.
\end{split}
\end{equation*}
Thus we have
\begin{equation}
Cb\dim{\frak m}=-2b\cdot
c({\frak k},\mathrm{ad}_{\frak p},\langle{\ },{\ }\rangle_{\frak k})
=1.
\end{equation}
Since
\begin{equation*}
b= \frac{\dim{\frak k}}
{\dim{\frak k}-(\dim{\frak p})\cdot
c({\frak k},\mathrm{ad}_{\frak p},\langle{\ },{\ }\rangle_{\frak k})}
\end{equation*}
by \cite[p.~591, Proposition 2.2]{Amar-Ohn03}, we obtain
\begin{equation*}
\begin{split}
b&=1-\frac{1}{2}\frac{\dim{\frak p}}{\dim{\frak k}},\\
-c({\frak k},\mathrm{ad}_{\frak p},
\langle{\ },{\ }\rangle_{\frak k})
&=\frac{\dim{\frak k}}{2\dim{\frak k}-\dim{\frak p}},\\
\Vert{\gamma_{1}}\Vert_{\frak u}^{2}&=
Cb=\frac{1}{\dim{\frak m}},\\
C&=\frac{2\dim{\frak k}}{(\dim{\frak m})(2\dim{\frak k}-\dim{\frak
p})}.
\end{split}
\end{equation*}
\end{proof}

\bigskip
\noindent
(1) The case $(U,K)=(SU(3),SO(3))$$\colon$

Define a basis $\{E_{1},E_{2},E_{3}\}$ of Lie algebra $\frak{su}(2)$
of $SU(2)$ by
\begin{equation*}
E_{1}:= \left(
\begin{array}{cc}
\sqrt{-1}&0\\
0&-\sqrt{-1}
\end{array}
\right),{\ }E_{2}:= \left(
\begin{array}{cc}
0&1\\
-1&0
\end{array}
\right),{\ }E_{3}:= \left(
\begin{array}{cc}
0&\sqrt{-1}\\
\sqrt{-1}&0
\end{array}
\right).
\end{equation*}
Let $\psi:SU(2)\rightarrow SO(3)$ be a universal covering Lie group
homomorphism defined by
\begin{equation*}
\mathrm{Ad}(a)(E_{1},E_{2},E_{3})=(E_{1},E_{2},E_{3})\psi(a).
\end{equation*}

Set $\tilde{K}:=SU(2)$, $\tilde{K}_{[{\frak
a}]}:=\psi^{-1}(K_{[{\frak a}]})$ and
$\tilde{K}_{0}:=\psi^{-1}(K_{0})$. Then we have
$\tilde{K}/\tilde{K}_{[{\frak a}]}\cong K/K_{[{\frak a}]}$ and
$\tilde{K}/\tilde{K}_{0}\cong K/K_{0}$. Explicitly the group $K_{0}$
is a finite subgroup of order $4$ generated by
\begin{equation*}
\left(
\begin{array}{ccc}
1&0&0\\
0&1&0\\
0&0&1\\
\end{array}
\right), \left(
\begin{array}{ccc}
1&0&0\\
0&-1&0\\
0&0&-1\\
\end{array}
\right), \left(
\begin{array}{ccc}
-1&0&0\\
0&-1&0\\
0&0&1\\
\end{array}
\right), \left(
\begin{array}{ccc}
-1&0&0\\
0&1&0\\
0&0&-1\\
\end{array}
\right)
\end{equation*}
and $K_{[{\frak a}]}$ is a finite subgroup of order $12$ generated
by $K_{0}$ and
\begin{equation*}
\left(
\begin{array}{ccc}
0&0&1\\
1&0&0\\
0&1&0\\
\end{array}
\right), \left(
\begin{array}{ccc}
0&1&0\\
0&0&1\\
1&0&0\\
\end{array}
\right) = \left(
\begin{array}{ccc}
0&0&1\\
1&0&0\\
0&1&0\\
\end{array}
\right)^{2}.
\end{equation*}

The group $\tilde{K}_{0}$ is a finite subgroup of order $8$
generated by
\begin{equation*}
\pm\left(
\begin{array}{cc}
1&0\\
0&1
\end{array}
\right), \pm\left(
\begin{array}{cc}
0&-1\\
1&0
\end{array}
\right), \pm\left(
\begin{array}{cc}
\sqrt{-1}&0\\
0&-\sqrt{-1}
\end{array}
\right), \pm\left(
\begin{array}{cc}
0&\sqrt{-1}\\
\sqrt{-1}&0
\end{array}
\right),
\end{equation*}
The group $\tilde{K}_{[{\frak a}]}$ is a finite subgroup of order
$24$ generated by $\tilde{K}_{[{\frak a}]}$ and
\begin{equation*}
\pm \left(
\begin{array}{cc}
\frac{1+\sqrt{-1}}{2}&\frac{1+\sqrt{-1}}{2}\\
\frac{-1+\sqrt{-1}}{2}&\frac{1-\sqrt{-1}}{2}
\end{array}
\right), \pm \left(
\begin{array}{cc}
\frac{-1+\sqrt{-1}}{2}&\frac{1+\sqrt{-1}}{2}\\
\frac{-1+\sqrt{-1}}{2}&\frac{-1-\sqrt{-1}}{2}
\end{array}
\right) = \pm \left(
\begin{array}{cc}
\frac{1+\sqrt{-1}}{2}&\frac{1+\sqrt{-1}}{2}\\
\frac{-1+\sqrt{-1}}{2}&\frac{1-\sqrt{-1}}{2}
\end{array}
\right)^{2}.
\end{equation*}

We know that
\begin{equation*}
{\mathcal D}(SU(2)) =\{(V_{m},\rho_{m}){\ }\vert{\ }m\in{\bold
Z},m\geq{0}\}.
\end{equation*}
Here $V_{m}$ denotes the complex vector space of complex homogeneous
polynomials of degree $m$ with two variables $z_{0},z_{1}$ and the
representation $\rho_{m}$ of $SU(2)$ on $V_{m}$ is defined by
\begin{equation*}
\left(\rho_{n}\left(
\begin{array}{cc}
a&-\bar{b}\\
b&\bar{a}
\end{array}
\right)f\right)(z_{0},z_{1}) =f((z_{0},z_{1}) \left(
\begin{array}{cc}
a&-\bar{b}\\
b&\bar{a}
\end{array}
\right))
\quad\text{ for each }
\left(
\begin{array}{cc}
a&-\bar{b}\\
b&\bar{a}
\end{array}
\right)\in{SU(2)}.
\end{equation*}
Set
\begin{equation*}
v^{(m)}_{k}:=\frac{1}{\sqrt{k!(m-k)!}}z_{0}^{m-k}z_{1}^{k}
\in{V_{m}} \quad(k=0,1,\dots,m)
\end{equation*}
and the standard Hermitian inner product $\langle\langle{\ },{\
}\rangle\rangle$ of $V_{m}$ invariant under $\rho_{m}$ is defined
such that $\{v^{(m)}_{0},\dots,v^{(m)}_{m}\}$ is a unitary basis of
$V_{m}$. Then there is a bijective correspondence between ${\mathcal
D}(SU(2))$ and $D(SU(2)) =\{m\Lambda_{1}{\ }\vert{\ }m\in{\bold
Z},m\geq{0}\}$, where $\Lambda_{1}$ denotes the fundamental weight
of $\frak{su}(2)$.

The Killing-Cartan form is given by
\begin{equation*}
B_{\frak{su}(2)}(X,Y)=4\mathrm{tr}(XY) \text{ for each
}X,Y\in{\frak{su}(2)}.
\end{equation*}

Let $\{X_{1},X_{2},X_{3}\}$ be an orthonormal basis of $\tilde{\frak
k}=\frak{su}(2)$ with respect to $\langle{\ },{\
}\rangle_{\tilde{\frak k}}$ defined by
\begin{equation*}
X_{1}=\frac{1}{2\sqrt{2}}E_{1},\quad
X_{2}=\frac{1}{2\sqrt{2}}E_{2},\quad X_{3}=\frac{1}{2\sqrt{2}}E_{3}.
\end{equation*}
The Casimir operator ${\mathcal C}$ of $SU(2)$ with respect to the
inner product $\langle{\ },{\ }\rangle_{\tilde{\frak k}}$ is given
as
\begin{equation*}
{\mathcal C}=\sum^{3}_{i=1}(X_{i})^{2}
=\frac{1}{8}\sum^{3}_{i=1}(E_{i})^{2}.
\end{equation*}
Then we obtain

\begin{lemma}[cf.~\cite{MOU84}]
\begin{enumerate}
\renewcommand{\labelenumi}{(\arabic{enumi})}
\item
The eigenvalue formula for the Casimir operator ${\mathcal C}$ is
given as
\begin{equation*}
\rho_{m}({\mathcal C})v=-\frac{m(m+2)}{8}v
\end{equation*}
for each $v\in{V_{m}}$.
\item
${\mathcal D}(\tilde{K},\tilde{K}_{0})$ is determined as follows :
$(\rho_{m},V_{m})\in{\mathcal D}(\tilde{K},\tilde{K}_{0})$ if and
only if $m$ is even and $m\geq{4}$. Then for each
$(\rho_{m},V_{m})\in{\mathcal D}(\tilde{K},\tilde{K}_{0})$ with
$m=2p$ for some integer $p$, the vector subspace
$(V_{m})_{\tilde{K}_{0}}$ is spanned by
\begin{equation*}
w_{i}:=\frac{1}{2}(v^{(m)}_{2(i-1)}+v^{(m)}_{4l-2(i-1)}){\ }
(i=1,\cdots,\ell+1)\quad \text{ if } p=2\ell,
\end{equation*}
or
\begin{equation*}
w^{\prime}_{i}:=\frac{1}{2}(v^{(m)}_{2i-1}-v^{(m)}_{4\ell-2i+3)}){\
} (i=1,\cdots,\ell)\quad \text{ if } p=2\ell+1.
\end{equation*}
\end{enumerate}
\end{lemma}

We must examine eigenvalues of ${\mathcal C}$ smaller than or equal
to $6$ for $(\rho_{m},V_{m})\in {\mathcal
D}(\tilde{K},\tilde{K}_{[{\frak a}]})\subset {\mathcal
D}(\tilde{K},\tilde{K}_{0})$. We observe that all the eigenvalues of
${\mathcal C}$ smaller than or equal to $6$ for $(\rho_{m},V_{m})\in
{\mathcal D}(\tilde{K},\tilde{K}_{0})$ are $3{\ }(m=4)$ and $6{\
}(m=6)$.

For $m=4$, we shall show that $(\rho_{4},V_{4})\notin{\mathcal
D}(\tilde{K},\tilde{K}_{[{\frak a}]}) {\ }(m=4)$. If we set
\begin{equation*}
B= \left(
\begin{array}{cc}
\frac{1+\sqrt{-1}}{2}&\frac{1+\sqrt{-1}}{2}\\
\frac{-1+\sqrt{-1}}{2}&\frac{1-\sqrt{-1}}{2}
\end{array}
\right) \in\tilde{K}_{[{\frak a}]},
\end{equation*}
then $B\notin\tilde{K}_{0}$ and
\begin{equation*}
(\rho_{4}(B))w_{1}
=-\frac{1}{4\sqrt{4!}}(z_{0}^{4}-6z_{0}^{2}z_{1}^{2}+z_{1}^{4})
\end{equation*}
\begin{equation*}
(\rho_{4}(B))w_{2}
=-\frac{1}{4\cdot{2!}}(z_{0}^{4}+2z_{0}^{2}z_{1}^{2}+z_{1}^{4}).
\end{equation*}
Thus if we assume that
$w=aw_{1}+bw_{2}\in(V_{m})_{\tilde{K}_{[{\frak a}]}}$, then by the
above equations $(\rho_{4}(B))w=w$ implies that $w=0$. Hence we
obtain $(V_{4})_{\tilde{K}_{[{\frak a}]}}=\{0\}$.

For $m=6$, since a simple computation implies that
$(\rho_{6}(B))w^{\prime}_{1}=\frac{1}{2}w^{\prime}_{1}$, we obtain
$(V_{6})_{\tilde{K}_{[{\frak a}]}} =\mathrm{span}_{\bold
C}\{w^{\prime}_{1}\}$ and thus $\dim{(V_{6})_{\tilde{K}_{[{\frak
a}]}}}=1$.

We conclude that $L={\mathcal G}(N^{3})$ is Hamiltonian stable and
the nullity is equal to $7=\dim(SO(5))-3$, and hence
$L={\mathcal G}(SO(3)/{\bold Z}_{2}+{\bold Z}_{2})$
is strictly Hamiltonian stable.

\bigskip
\noindent
(2) The case $(U,K)=(SU(3)\times SU(3),SU(3))$$\colon$
Then
\begin{equation*}
\begin{split}
&K=SU(3),{\ }K_{0}=T^{2} \\
&D(K,K_{0})=D(SU(3),T^{2}),{\ } D(K,K_{[\frak
a]})=D(SU(3),T^{2}\cdot{\bold Z}_{3})
\end{split}
\end{equation*}
and
\begin{equation*}
D(SU(3),T^{2}\cdot{\bold Z}_{3})\subset{D(SU(3),T^{2})}.
\end{equation*}
Let $\{\alpha_{1},\alpha_{2}\}$ be the fundamental root system of
$SU(3)$ and $\{\Lambda_{1},\Lambda_{2}\}$ be the fundamental weight
system of $SU(3)$. We use the results of Satoru Yamaguchi
\cite{SYamaguchi79} as follows : Each $\Lambda\in{D(SU(3),T^{2})}$
can be uniquely expressed as
\begin{equation*}
\Lambda=\sum_{i}m_{i}\Lambda_{i}=\sum_{i}p_{i}\alpha_{i},
\end{equation*}
where $m_{i}\in{\bold Z}$, $m_{i}\geq{0}$, $p_{i}\in{\bold Z}$,
$p_{i}\geq{1}$ and
\begin{equation*}
m_{1}=2p_{1}-p_{2}\geq{0},\quad m_{2}=-p_{1}+2p_{2}\geq{0}.
\end{equation*}
The eigenvalue formula is
\begin{equation*}
-c(\Lambda,\langle{\ },{\ }\rangle_{\frak k})
=\frac{1}{6}\left(m_{1}p_{1}+m_{2}p_{2}+2p_{1}+2p_{2}\right)
\end{equation*}
for each $\Lambda\in{D(SU(3),T^{2})}$. Therefore we get
\begin{equation*}
\begin{split}
&\{\Lambda\in{D(SU(3),T^{2})}{\ }\vert{\ }
-c(\Lambda,\langle{\ },{\ }\rangle_{\frak k})\leq{2}\}\\
=&\{{\ }0,{\ }3\Lambda_{1}{\ }((p_{1},p_{2})=(2,1)),
{\ }3\Lambda_{2}{\ }((p_{1},p_{2})=(1,2)),\\
&{\ }\Lambda_{1}+\Lambda_{2}{\ }((p_{1},p_{2})=(1,1)){\ }\}.
\end{split}
\end{equation*}
Since
$\Lambda_{1}+\Lambda_{2}\notin{D(SU(3),T^{2}\cdot{\bold Z}_{3})}$,
we obtain that
\begin{equation*}
\begin{split}
&\{\Lambda\in{D(SU(3),T^{2}\cdot{\bold Z}_{3})}{\ }\vert{\ }
-c(\Lambda,\langle{\ },{\ }\rangle_{\frak k})\leq{2}\}\\
=&\{{\ }0,{\ }3\Lambda_{1}{\ }((p_{1},p_{2})=(2,1)), {\
}3\Lambda_{2}{\ }((p_{1},p_{2})=(1,2)){\ }\}.
\end{split}
\end{equation*}
and $-c(3\Lambda_{1},\langle{\ },{\ }\rangle_{\frak k})
=-c(3\Lambda_{2},\langle{\ },{\ }\rangle_{\frak k}) =2$.
By \cite{McKay-Patera} we have
$\dim_{\bold C}(V_{3\Lambda_{1}})_{T^{2}}=
\dim_{\bold C}(V_{3\Lambda_{2}})_{T^{2}}=1$.
Hence
the nullity is equal to
$\dim(V_{3\Lambda_{1}})+\dim(V_{3\Lambda_{2}})=10+10=20
=\dim(SO(8))-\dim(SU(3))$. We conclude that
$L={\mathcal G}(SU(3)/T^{2})$ is strictly Hamiltonian stable.

\bigskip
\noindent
(3) The case $(U,K)=(SU(6),Sp(3))$$\colon$
Then
\begin{equation*}
\begin{split}
&K=Sp(3),{\ }K_{0}=Sp(1)^{3},\\
&D(K,K_{0})=D(Sp(3),Sp(1)^{3}),{\ } D(K,K_{[\frak
a]})=D(Sp(3),Sp(1)^{3}\cdot{\bold Z}_{3})
\end{split}
\end{equation*}
and
\begin{equation*}
D(Sp(3),Sp(1)^{3}\cdot{\bold Z}_{3})\subset{D(Sp(3),Sp(1)^{3})}
\subset{D(Sp(3),T^{3})}.
\end{equation*}
Let $\{\alpha_{1},\alpha_{2},\alpha_{3}\}$ be the fundamental root
system of $Sp(3)$ and $\{\Lambda_{1},\Lambda_{2},\Lambda_{3}\}$ be
the fundamental weight system of $Sp(3)$. The results from
\cite{SYamaguchi79} are as follows : Each
$\Lambda\in{D(Sp(3),T^{3})}$ can be uniquely expressed as
\begin{equation*}
\Lambda=\sum_{i}m_{i}\Lambda_{i}=\sum_{i}p_{i}\alpha_{i},
\end{equation*}
where $m_{i}\in{\bold Z}$, $m_{i}\geq{0}$, $p_{i}\in{\bold Z}$,
$p_{i}\geq{1}$ and
\begin{equation*}
m_{1}=2p_{1}-p_{2}\geq{0},{\ } m_{2}=-p_{1}+2p_{2}-p_{3}\geq{0},{\ }
m_{3}=-p_{2}+2p_{3}\geq{0}.
\end{equation*}
The eigenvalue formula is
\begin{equation*}
-c(\Lambda,\langle{\ },{\ }\rangle_{\frak k}) =\frac{1}{16}
\left(m_{1}p_{1}+m_{2}p_{2}+2m_{3}p_{3}+ 2p_{1}+2p_{2}+4p_{3}\right)
\end{equation*}
for each $\Lambda\in{D(Sp(3),T^{3})}$. Therefore we get
\begin{equation*}
\begin{split}
&\{{\ }\Lambda\in{D(Sp(3),T^{3})}{\ }\vert{\ }
-c(\Lambda,\langle{\ },{\ }\rangle_{\frak k})\leq{3/2}{\ }\}\\
=&\{{\ }0,{\ }2\Lambda_{1}{\ }((p_{1},p_{2},p_{3})=(2,2,1)),\\
&{\ }\Lambda_{2}{\ }((p_{1},p_{2},p_{3})=(1,2,1)), {\
}\Lambda_{1}+\Lambda_{3}{\ }((p_{1},p_{2},p_{3})=(2,3,2)){\ }\}.
\end{split}
\end{equation*}
Since we see that $2\Lambda_{1}\notin{D(Sp(3),Sp(1)^{3})}$ by using
\cite{McKay-Patera}, we obtain that
\begin{equation*}
\{{\ }\Lambda\in{D(Sp(3),Sp(1)^{3}\cdot{\bold Z}_{3})}{\ }\vert{\ }
-c(\Lambda,\langle{\ },{\ }\rangle_{\frak k})\leq{3/2}{\ }\} =\{{\
}0,{\ }\Lambda_{1}+\Lambda_{3}{\ }\}.
\end{equation*}
and
$-c(\Lambda_{1}+\Lambda_{3},\langle{\ },{\ }\rangle_{\frak k})={3/2}$.
By \cite{McKay-Patera} we have
$\dim_{\bold C}(V_{\Lambda_{1}+\Lambda_{3}})_{Sp(1)^{3}}=1$.
Hence
the nullity is equal to
$\dim(\Lambda_{1}+\Lambda_{3})=70
=\dim(SO(14))-\dim(Sp(3))(=91-21)$. We conclude that $L={\mathcal
G}(Sp(3)/Sp(1)^{3})$ is strictly Hamiltonian stable.

\bigskip
\noindent
(4) The case $(U,K)=(E_{6},F_{4})$$\colon$
Then
\begin{equation*}
\begin{split}
&K=F_{4},{\ }K_{0}=Spin(8), \\
&D(K,K_{0})=D(F_{4},Spin(8)),{\ } D(K,K_{[\frak
a]})=D(F_{4},Spin(8)\cdot{\bold Z}_{3})
\end{split}
\end{equation*}
and
\begin{equation*}
D(F_{4},Spin(8)\cdot{\bold Z}_{3})\subset{D(F_{4},Spin(8))}
\subset{D(F_{4},T^{4})}.
\end{equation*}
Let $\{\alpha_{1},\alpha_{2},\alpha_{3},\alpha_{4}\}$ be the
fundamental root system of $F_{4}$ and
$\{\Lambda_{1},\Lambda_{2},\Lambda_{3},\Lambda_{4}\}$ be the
fundamental weight system of $F_{4}$. The results from
\cite{SYamaguchi79} are as follows : Each
$\Lambda\in{D(F_{4},T^{4})}$ can be uniquely expressed as
\begin{equation*}
\Lambda=\sum_{i}m_{i}\Lambda_{i}=\sum_{i}p_{i}\alpha_{i},
\end{equation*}
where $m_{i}\in{\bold Z}$, $m_{i}\geq{0}$, $p_{i}\in{\bold Z}$,
$p_{i}\geq{1}$ and
\begin{equation*}
\begin{split}
&m_{1}=2p_{1}-p_{2}\geq{0},{\ }
m_{2}=-p_{1}+2p_{2}-p_{3}\geq{0},\\
&m_{3}=-2p_{2}+2p_{3}-p_{4}\geq{0},{\ } m_{4}=-p_{3}+2p_{4}\geq{0}.
\end{split}
\end{equation*}
The eigenvalue formula is
\begin{equation*}
-c(\Lambda,\langle{\ },{\ }\rangle_{\frak k}) =\frac{1}{18}
\left(m_{1}p_{1}+m_{2}p_{2}+\frac{1}{2}m_{3}p_{3}+\frac{1}{2}m_{4}p_{4}
+2p_{1}+2p_{2}+p_{3}+p_{4}\right)
\end{equation*}
for each $\Lambda\in{D(F_{4},T^{4})}$. Therefore we get
\begin{equation*}
\begin{split}
&\{{\ }\Lambda\in{D(F_{4},T^{4})}{\ }\vert{\ }
-c(\Lambda,\langle{\ },{\ }\rangle_{\frak k})\leq{4/3}{\ }\}\\
=&\{{\ }0,
{\ }\Lambda_{1}{\ }((p_{1},p_{2},p_{3},p_{4})=(2,3,4,2)),\\
&{\ }\Lambda_{3}{\ }((p_{1},p_{2},p_{3},p_{4})=(2,4,6,3)), {\
}\Lambda_{4}{\ }((p_{1},p_{2},p_{3},p_{4})=(1,2,3,2)){\ }\}.
\end{split}
\end{equation*}
Since $\Lambda_{4}\notin{D(F_{4},Spin(8)\cdot{\bold Z}_{3})}$ and we
see $\Lambda_{1}\notin{D(F_{4},Spin(8))}$ by using
\cite{McKay-Patera} we obtain that
\begin{equation*}
\begin{split}
&\{{\ }\Lambda\in{D(F_{4},Spin(8)\cdot{\bold Z}_{3})}{\ }\vert{\ }
-c(\Lambda,\langle{\ },{\ }\rangle_{\frak k})\leq{4/3}{\ }\}\\
=&\{{\ }0,{\ }\Lambda_{3}{\ }((p_{1},p_{2},p_{3},p_{4})=(2,4,6,3)){\
}\}.
\end{split}
\end{equation*}
and $-c(\Lambda_{3},\langle{\ },{\ }\rangle_{\frak k})=4/3$.
By \cite{McKay-Patera} we have
$\dim_{\bold C}(V_{\Lambda_{3}})_{Spin(8)}=1$.
Hence
the nullity is equal to $\dim(\Lambda_{3})=273
=\dim(SO(26))-\dim(F_{4})(=325-52)$. We conclude that $L={\mathcal
G}(F_{4}/Spin(8))$ is strictly Hamiltonian stable. This proves the
theorem.


\end{document}